\documentclass[10pt,a4paper]{article}
\usepackage[utf8]{inputenc}
\usepackage{amsmath}
\usepackage{amsfonts}
\usepackage{amssymb}
\usepackage{amsthm}
\usepackage{setspace}
\usepackage{graphicx}
\usepackage{hyperref}
\usepackage{mathtools}
\usepackage{subcaption}
\usepackage[T1]{fontenc}
\usepackage{xcolor}
\usepackage{float}
\usepackage[noadjust]{cite}
\usepackage{cite}
\usepackage[margin=1.2in]{geometry}

\graphicspath{ {./images/} }
\newtheorem{thm}{Theorem}[section]
\newtheorem{prop}[thm]{Proposition}
\newtheorem{lemma}[thm]{Lemma}

\newtheorem{algo}[thm]{Algorithm}
\newtheorem*{theorem*}{Theorem A}
\newtheorem*{theorem**}{Theorem D}
\newtheorem*{coro*}{Corollary B}
\newtheorem*{coro**}{Corollary C}
\newtheorem*{prop*}{Proposition E}

\theoremstyle{definition}
\newtheorem{defi}[thm]{Definition}
\newtheorem{ex}[thm]{Example}
\newtheorem{rem}[thm]{Remark}
\newtheorem{conj}[thm]{Conjecture}

\newcommand{\quotient}[2]{{\raisebox{.2em}{$#1$}\left/\raisebox{-.2em}{$#2$}\right.}}

\begin{document}

\begin{center}
{\large \textbf{Acylindrical hyperbolicity for Artin groups of dimension 2}}
\end{center}

\begin{center}
Nicolas Vaskou
\footnote{\text{E-mail address:} \texttt{\href{mailto: ncv1@hw.ac.uk}{ncv1@hw.ac.uk}}

2020 \textit{Mathematics subject classification.} 20F65, 20F36, 20F67.

\textit{Key words.} Artin groups, acylindrically hyperbolic groups, CAT(0) geometry, Deligne complex.}
\end{center}

\begin{abstract}
\noindent In this paper, we show that every irreducible $2$-dimensional Artin group $A_{\Gamma}$ of rank at least $3$ is acylindrically hyperbolic. We do this by studying the action of $A_\Gamma$ on its modified Deligne complex. Along the way, we prove results of independent interests on the geometry of links of this complex.
\end{abstract}

\section{Introduction}

Artin-Tits groups were introduced by Tits (\cite{tits1966normalisateurs}), as “extended Coxeter groups”. Consider a simplicial graph $\Gamma$ with vertex set $V(\Gamma)$, and suppose that every edge $e^{ab} \in E(\Gamma)$ between two vertices $a$ and $b$ has integer coefficient $m_{ab} \geq 2$. The Artin group defined by $\Gamma$ is the group with presentation
$$A_{\Gamma} \coloneqq \langle \ V(\Gamma) \ | \ \underbrace{aba\cdots}_{m_{ab} \text{ terms}} = \underbrace{bab\cdots}_{m_{ab} \text{ terms}} \text{ for every edge } e^{ab} \in E(\Gamma) \ \rangle. \ (*)$$
It is common to say that $m_{ab} = \infty$ if $a$ and $b$ are not connected by an edge. The rank of $A_{\Gamma}$ is the cardinality of $V(\Gamma)$, and is assumed to be finite. Similarly, there is a Coxeter group $W_{\Gamma}$ attached to $\Gamma$, whose presentation is obtained from $(*)$ by adding $s^2 = 1$ for every generator $s$ in $V(\Gamma)$. The two groups are connected through the natural projection $A_{\Gamma} \twoheadrightarrow W_{\Gamma}$ that restricts to the identity on $S$. The class of Artin groups encompasses a large spectrum of groups, going from free groups to free abelian groups, and including braid groups (the braid relation corresponds to a coefficient $3$).

Coxeter groups are well understood. For instance, they are virtually torsion-free (\cite{davis2008geometry}), have solvable word and conjugacy problem (this is a consequence of them being CAT(0) groups (\cite{moussong1988hyperbolic})), and their centre are finite (isomorphic to $(\mathbf{Z}/2\mathbf{Z})^n$ for some $n \geq 0$) (\cite{hosaka2005center}). On the other hand, most results obtained for Artin groups only concern certain classes. For instance, the following conjectures remain open, at least in general (see \cite{charney2016problems} for more questions about Artin groups):

\begin{conj} \label{ConjAlg} Consider an Artin group $A_{\Gamma}$. Then:
\\(1) $A_{\Gamma}$ has solvable word and conjugacy problem.
\\(2) $A_{\Gamma}$ is torsion-free.
\\(3) The centre of $A_{\Gamma}$ is trivial (if $A_{\Gamma}$ is non-spherical) or infinite cyclic (if $A_{\Gamma}$ is spherical).
\\(4) $A_{\Gamma}$ satisfies the $K(\pi,1)$ conjecture.
\end{conj}

\noindent We give here three classes of Artin groups for which substantial progress has been made:
\\$\bullet$ Spherical Artin groups, namely Artin groups $A_{\Gamma}$ whose associated Coxeter group $W_{\Gamma}$ is finite.
\\$\bullet$ Artin groups of type FC, namely Artin groups $A_{\Gamma}$ for which every complete subgraph $\Gamma' \subseteq \Gamma$ generates an Artin subgroup $A_{\Gamma'}$ of spherical type. Note that spherical Artin groups and right angled Artin groups (those whose coefficients are either $2$ or $\infty$) are of type FC.
\\$\bullet$ Artin groups of dimension $2$, where the dimension of an Artin group is the maximal rank of a spherical Artin subgroup. This class includes the class of large Artin groups (those with coefficients at least $3$). Artin groups of dimension $2$ are known to be exactly those (see \cite{charney1995k}) for which ($\Gamma$ is not discrete and) every triangle in $\Gamma$ with vertices $a, b, c$ satisfies
$$\frac{1}{m_{ab}} + \frac{1}{m_{ac}} + \frac{1}{m_{bc}} \leq 1.$$
Spherical Artin groups are well understood. The existence of a normal form (\cite{garside1969braid},\cite{elrifai1994algorithms},\cite{dehornoy1999gaussian}) and of a contractible simplicial complex on which the group acts in a good manner (referred as the Deligne complex and due to Deligne (\cite{deligne1972immeubles})) has highly helped to understand their structure and geometry. In particular, they are known to satisfy all the above conjectures. Artin groups of type FC and of dimension $2$ are well understood and satisfy all of the above conjectures as well (\cite{charney1995k}, \cite{chalopin2020helly}, \cite{huang2019metric}).

In more recent years, there has been an increasing interest in understanding the geometry of Artin groups. In particular, the action of an Artin group on a space with properties that encodes some kind of non-positive curvature, such as being CAT(0), systolic or acylindrically hyperbolic, turns out to be particularly interesting. However, the geometry of Artin groups remains very mysterious. For instance, it is not known in general whether the following conjectures hold:

\begin{conj} \label{ConjGeom} Consider an Artin group $A_{\Gamma}$. Then:
\\(1) $A_{\Gamma}$ is CAT(0), i.e. acts properly and cocompactly on a CAT(0) space.
\\(2) The central quotient $A_{\Gamma} / Z(A_{\Gamma})$ is acylindrically hyperbolic.
\end{conj}

As an example, Coxeter groups are known to be CAT(0) (\cite{moussong1988hyperbolic}). There are partial results to the above conjectures. For instance, Conjecture \ref{ConjGeom}.(1) is known to hold for right angled Artin groups (\cite{charney1995k}), some classes of $2$-dimensional Artin groups (\cite{brady2002two},\cite{brady2000three},\cite{haettel2019xxl}) or spherical Artin groups of rank $3$ (\cite{brady2000artin}). As regards to Conjecture \ref{ConjGeom}.(2), it is known to hold for Artin groups of spherical type (\cite{calvez2017acylindrical}). It is then enough to look at what happens when the group is non-spherical. In that case, the centre $Z(A_{\Gamma})$ is conjectured to be trivial (Conjecture \ref{ConjAlg}.(3)), and hence the question essentially comes down to asking whether $A_{\Gamma}$ is acylindrically hyperbolic.

The notion of acylindricity goes back to (\cite{sela1997acylindrical}) and gives conditions on the size of stabilisers associated with group actions on trees. In the more general case of metric spaces, the definition is due to (\cite{bowditch2008tight}): a group $G$ is said to act acylindrically on a space $X$ if for every $R \geq 0$, there exist $N > 0$, $L > 0$ such that
$$\forall x, y \in X, \ d(x,y) \geq L \Rightarrow |\{g \in G \ | \ d(x,gx) \leq R, \ d(y,gy) \leq R \}| \leq N.$$
A group $G$ is said to be acylindrically hyperbolic if it is not virtually cyclic and has an acylindrical action on a hyperbolic space. The condition of acylindrical hyperbolicity introduced by Osin (\cite{osin2016acylindrically}) merges many previously known results, bringing together classes such as mapping class groups, Out($F_n$) for $n \geq 2$, many CAT(0) groups and most of $3$-manifold groups. Nevertheless, acylindrical hyperbolicity is still strong enough to ensure interesting properties for the group: for instance, acylindrically hyperbolic groups are SQ-universal and contains free normal subgroups.

Many classes of Artin groups are known to be acylindrically hyperbolic. For instance, right angled Artin groups that are not cyclic nor reducible are acylindrically hyperbolic (\cite{osin2016acylindrically}). In \cite{charney2019artin}, Charney and Morris-Wright also showed that Artin groups $A_{\Gamma}$ whose defining graph $\Gamma$ is not a join are acylindrically hyperbolic, extending the result of Chatterji and Martin for Artin groups of type FC whose defining graph has diameter at most $3$ (\cite{chatterji2016note}). It is also known from \cite{haettel2019xxl} that XXL Artin groups (those with coefficient at least $5$) of rank at least $3$ are acylindrically hyperbolic. More recently, Kato and Oguni also showed that triangle-free Artin groups and Artin groups of large type associated to cones over square-free bipartite graphs are also acylindrically hyperbolic (\cite{kato2020acylindrical}). It was also proved in \cite{calvez2020euclidean} that Artin groups of euclidean type are acylindrically hyperbolic. In \cite{martin2019acylindrical}, Martin and Przytycki showed that $2$-dimensional Artin groups of hyperbolic type (whose associated Coxeter groups are hyperbolic) are acylindrically hyperbolic.

Let us recall that an Artin group $A_{\Gamma}$ is said to be reducible if $\Gamma$ is a join of two non-trivial subgraphs such that every edge of the join has coefficient $2$. Note that if $A_{\Gamma}$ is reducible, then it can be written as a direct product $A_{\Gamma_1} \times A_{\Gamma_2}$ of the groups generated by the two subgraphs, hence is not acylindrically hyperbolic (\cite{osin2016acylindrically}). The goal of the present paper is to show the following extension to some of the previously mentioned results, answering to the question of acylindrical hyperbolicity for Artin groups of dimension $2$ in general:

\begin{theorem*} Every irreducible $2$-dimensional Artin group of rank at least $3$ is acylindrically hyperbolic.
\end{theorem*}

For instance, it was not known whether the rather simple following Artin group was acylindrically hyperbolic:

\begin{figure}[H]
\centering
\includegraphics[scale=0.5]{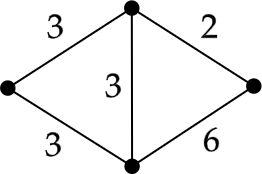}
\caption{Example of an Artin group for which acylindrical hyperbolicity was not previously known.}
\label{FigExOfAH}
\end{figure}

Because direct products of infinite groups are never acylindrically hyperbolic, it is clear from the previous Theorem that the following holds:

\begin{coro*} Irreducible $2$-dimensional Artin groups of rank at least $3$ cannot be decomposed as direct products of infinite groups.
\end{coro*}

Recall that $2$-dimensional Artin groups are torsion-free. Moreover, acylindrically hyperbolic groups have finite centres (Corollary 7.2 of \cite{osin2016acylindrically}). Thus it follows from Theorem A that irreducible $2$-dimensional Artin groups of rank at least $3$ have trivial centres. This also holds if $A_{\Gamma}$ is reducible, as $2$-dimensional reducible Artin groups are direct products of free groups, hence have trivial centres. If $A_{\Gamma}$ is irreducible and has rank $2$ then it is a dihedral Artin group with coefficient at least $3$, and $A_{\Gamma} / Z(A_{\Gamma})$ is virtually a free group (\cite{brady2000three},\cite{ciobanu2020equations}), hence acylindrically hyperbolic.

Putting together everything that we just discussed, we get the following corollary. In particular, we prove Conjecture \ref{ConjGeom}.(2) for $2$-dimensional Artin groups, and we give a new proof of Conjecture \ref{ConjAlg}.(3), which could already be deduced from \cite{godelle2007artin} although it is not explicitly stated:

\begin{coro**} Artin groups of dimension $2$ and rank at least $3$ have trivial centre. Irreducible Artin groups $A_{\Gamma}$ of dimension $2$ have acylindrically hyperbolic central quotient $A_{\Gamma} / Z(A_{\Gamma})$.
\end{coro**}

Checking whether an action is acylindrical can be tough, as it essentially comes down to controlling the geodesics between two metric balls. Instead, one usually looks for an “acylindrical direction” in the action, more precisely a WPD element with a strongly contracting orbit, from which one can construct an acylindrical action on a larger space (Theorem H,\cite{bestvina2015constructing}). That approach was followed to prove the acylindrical hyperbolicity of different classes of groups (\cite{bestvina3730hyperbolic},\cite{gruber2018infinitely},\cite{minasyan2015acylindrical}). That said, this condition remains hard to check when the space acted upon is not locally compact.

In \cite{charney1995k}, Charney and Davis constructed an analogue of the Deligne complex for non-spherical Artin groups, called the modified Deligne complex, which turns out to be equally interesting. It is conjectured that this space always supports a CAT(0) metric. In particular, it was proved that non-spherical Artin groups whose modified Deligne complex is CAT(0) satisfy Conjecture \ref{ConjAlg}.(2) and Conjecture \ref{ConjAlg}.(4). So far, it has been shown that this conjecture holds for Artin groups of type FC and Artin groups of dimension $2$ (\cite{charney1995k}).

In this paper, we focus on the action of Artin groups of dimension $2$ on their modified Deligne complexes. In general, this action is not acylindrical and the space is not hyperbolic. Unfortunately, the modified Deligne complex is not locally compact either, which make the use of the above mentioned WPD condition harder. To bypass that problem we will use a criterion from Martin (\cite{martin2017acylindrical}) that uses a variant of the WPD condition, generalizing to higher dimension a result of \cite{minasyan2015acylindrical} for groups acting on trees. This criterion  ensures that under specific hypotheses, the existence of weakly malnormal local groups associated with actions in non-positive curvature ensure that the group is acylindrically hyperbolic. We recall that criterion thereafter (in a slightly more specific case):

\begin{theorem**} \textbf{(\cite{martin2017acylindrical}, Theorem B)} Let $X$ be a CAT(0) simplicial complex, together with an action by simplicial isomorphisms of a group $G$. Assume that there exists a vertex $v$ of $X$ with stabiliser $G_v$ such that:
\\(1) The orbits of $G_v$ on the link $Lk_X(v)$ are unbounded, for the associated angular metric.
\\(2) $G_v$ is weakly malnormal in $G$, i.e. there exists an element $g \in G$ such that $G_v \cap g G_v g^{-1}$ is finite.
\\Then $G$ is either virtually cyclic or acylindrically hyperbolic.
\end{theorem**}

The proof of Theorem A has two major steps. First, we show that if $A_{\Gamma}$ is not right angled then there exists a vertex $v$ in the modified Deligne complex $D_{\Gamma}$ associated to the Artin group $A_{\Gamma}$ that satisfies Theorem D.(1). Then, we show geometrically that the stabiliser of this vertex is weakly malnormal in $A_{\Gamma}$, satisfying Theorem D.(2). The result then follows from Theorem D.

The paper is organised as follows. In section $2$, we recall some basic definitions, including complexes of groups and their associated developments. In section $3$, we talk about the modified Deligne complex and give a more precise description of its geometry, and that of its links. In section $4$, we show the following result:

\begin{prop*} Let $A_{\Gamma}$ be a $2$-dimensional Artin group of rank at least $3$ with modified Deligne complex $D_{\Gamma}$. Suppose that there exists a vertex $v_{ab} \in D_{\Gamma}$ whose stabiliser $A_{ab}$ has coefficient $3 \leq m_{ab} < \infty$. Then:
\\(1) The orbits of $A_{ab}$ on $Lk_{D_{\Gamma}}(v_{ab})$ are unbounded.
\\(2) More precisely, the orbits of $\langle g \rangle$ on $Lk_{D_{\Gamma}}(v_{ab})$ are quasi-isometrically embedded if and only if $g \in A_{ab}$ is not trivial, nor the conjugate of a power of a standard generator $a$ or $b$.
\end{prop*}

Note that when applied to our specific case, the first hypothesis of Theorem D is exactly the first result of Proposition E, so that we a priori don't need to prove Proposition E.(2). However, Proposition E.(2) remains interesting on its own, as it will be for instance used by Hagen-Martin-Sisto to prove that extra-large type Artin groups are virtually hierarchically hyperbolic (\cite{hagen2021extra}). In section $5$, we reduce the question of asking whether a dihedral subgroup of $A_{\Gamma}$ is weakly malnormal to a geometric question (see Lemma \ref{StabGeod}). The existence of weakly malnormal subgroups turns out to be implied by a simple geometric condition on the geodesics in the complex. We are able to show that this condition holds for all irreducible $2$-dimensional Artin groups of rank at least $3$ (assuming they are not free nor right angled, see Lemma \ref{LemmaGeod}). In particular, we show that the local group $G_v$ is weakly malnormal in $A_{\Gamma}$, i.e. that $G_v$ satisfies Theorem D.(2). We can then use Theorem D and prove Theorem A as an immediate consequence.
\bigskip

\textbf{Acknowledgments:} I would like to thank my supervisor Alexandre Martin with whom I had many constructive discussions. This work was partially supported by the EPSRC New Investigator Award EP/S010963/1. I would also like to thank the anonymous referee for their useful remarks and advices, and for proposing a strategy for improving and making optimal Proposition E.(2).

\section{Preliminaries}

In this section we recall general notions that we will use throughout the paper. In particular, we talk about complexes of groups and their associated developments. We begin by making a small remark about the definition of an Artin group that is given in the introduction:

\begin{rem} (1) Let $\Gamma$ be a graph defining an Artin group $A_{\Gamma}$, and let $\Gamma'$ be a subgraph of $\Gamma$. It is known from \cite{van1983homotopy} that the subgroup of $\Gamma$ generated by $V(\Gamma')$ is isomorphic to the Artin group $A_{\Gamma'}$, hence we will just write $A_{\Gamma'}$ to talk about the subgroup of $A_{\Gamma}$ generated by $V(\Gamma')$.
\\(2) Let $A_{\Gamma}$ be either a free group on two generators $a$ and $b$ or a dihedral Artin group, i.e. an Artin group defined by a graph $\Gamma$ that only has two vertices $a$ and $b$ and an edge between them. Then we will just write $A_{ab}$ instead of $A_{\Gamma}$.
\end{rem}

We now come to the definitions of complexes of groups and their associated developments. For more information about these notions, we refer the reader to (\cite{bridson2013metric}, Chapter II.12).

\begin{defi} A \textbf{simple complex of groups} $G(\mathcal{Q})$ over a poset (partially ordered set) $\mathcal{Q}$ consists of:
\\(1) For each element $\sigma \in \mathcal{Q}$, a group $G_{\sigma}$ called the local group at $\sigma$.
\\(2) For each $\tau < \sigma$, an injective morphism $\psi_{\tau \sigma} : G_{\sigma} \hookrightarrow G_{\tau}$ such that
$$\tau < \sigma < \rho \Longrightarrow \psi_{\tau \rho} = \psi_{\tau \sigma} \psi_{\sigma \rho},$$
i.e. every diagram of maps commute. A \textbf{simple morphism} $\varphi$ from $G(\mathcal{Q})$ to a group $G$ is a map that associates to each $\sigma \in \mathcal{Q}$ a morphism $\varphi_{\sigma} : G_{\sigma} \rightarrow G$ such that if $\tau < \sigma$ then $\varphi_{\sigma} = \varphi_{\tau} \psi_{\tau \sigma}$. The map $\varphi$ is said to be injective on the local groups if $\varphi_{\sigma}$ is injective for each $\sigma \in \mathcal{Q}$. The group
$$\widehat{G(\mathcal{Q})} \coloneqq \lim\limits_{\overset{\longrightarrow}{\sigma \in \mathcal{Q}}} G_{\sigma},$$
which is the direct limit of the system $(G_{\sigma},\psi_{\tau \sigma})$, is called the \textbf{fundamental group} of the complex of groups.
\end{defi}

\begin{ex} (1) A $n$-dimensional simplex of groups is a complex of groups over the poset of the faces of a simplex of dimension $n$. More precisely, if $\Delta$ is a simplex of dimension $n$ with faces $F_1,\cdots,F_n$, a face of codimension $k$ in $\Delta$ can be uniquely written as
$$F_I \coloneqq \bigcap\limits_{i \in I} F_i,$$
where $I \subseteq \{1,\cdots,n\}$ is such that $|I| = k$. Then we associate with every face $F_I$ a local group $G_{F_I}$ and to every inclusion $F_I \subseteq F_J$ an injective morphism  $\psi_{F_I F_J} : G_{F_J} \hookrightarrow G_{F_I}$. A triangle of groups is a simplex of groups with $n = 2$ (see Figure \ref{FigTriOfGrp} for an example).

(2) Let $A_{\Gamma}$ be an Artin group of rank $3$. The proper Artin subgroups $A_{\Gamma'}$ of $A_{\Gamma}$ form a poset that gives rise to a triangle of groups whose fundamental group is precisely $A_{\Gamma}$:
\begin{figure}[H]
\centering
\includegraphics[scale=1.1]{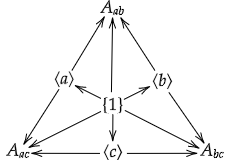}
\caption{Triangle of Artin groups with fundamental group $A_{\Gamma}$. The local group corresponding to the face is $\{1\}$, those corresponding to the edges are the cyclic groups $\langle a \rangle \cong \mathbf{Z}$, $\langle b \rangle \cong \mathbf{Z}$ and $\langle c \rangle \cong \mathbf{Z}$, and those corresponding to the vertices are the dihedral Artin groups $A_{ab}$, $A_{ac}$ and $A_{bc}$. The maps are just the natural inclusions.}
\label{FigTriOfGrp}
\end{figure}
\end{ex}

\begin{defi} \label{DefiDev} Let $X$ be a simplicial complex and let $\mathcal{P}$ be the set of its simplices, partially ordered by inclusion. Let $G$ be a group acting simplicially on $X$ with strict fundamental domain $Y$, and let $\mathcal{Q} \coloneqq \{\sigma \in \mathcal{P} \ | \ \sigma \subseteq Y \}$. From that we can recover a complex of groups $G(\mathcal{Q})$ in the following way :
\\$\bullet$ To each element $\sigma \in \mathcal{Q}$ corresponds a subgroup $G_{\sigma}$ of $G$ which is the stabiliser of $\sigma$ through the action of $G$.
\\$\bullet$ To every inclusion $\tau \subseteq \sigma \in \mathcal{Q}$ corresponds a map $\psi_{\tau \sigma} : G_{\sigma} \hookrightarrow G_{\tau}$ that is just the natural inclusion of the corresponding stabilisers.
\\The complex of groups $G(\mathcal{Q})$ is then defined to be
$$G(\mathcal{Q}) \coloneqq \{(G_{\sigma},\psi_{\tau \sigma}) \ | \ \sigma, \tau \in \mathcal{Q}, \ \tau \subseteq \sigma \}.$$
Notice that the inclusions $\varphi_{\sigma} : G_{\sigma} \rightarrow G$ give a simple morphism $\varphi : G(\mathcal{Q}) \rightarrow G$ that is injective on the local groups. A complex of groups $G(\mathcal{Q})$ is said to be \textbf{developable} if there exists a simplicial complex $X$ and a simplicial action of $G$ on $X$ with strict fundamental domain some subcomplex $Y$, such that the complex of groups recovered in the previous way is precisely $G(\mathcal{Q})$.
\end{defi}

\begin{defi} Let $Y$ be a simplicial complex and let $\mathcal{Q}$ be the set of its simplices, partially ordered by inclusion. Let now $G(\mathcal{Q})$ be a complex of groups and let $\varphi : G(\mathcal{Q}) \rightarrow G$ be a simple morphism to some group $G$, that is injective on the local groups. The \textbf{development} $D(Y,\varphi)$ of $Y$ along $\varphi$ is defined by
$$D(Y,\varphi) \coloneqq \quotient{G \times Y}{\sim},$$
where $(g,x) \sim (g',x') \Longleftrightarrow x = x'$ and $g^{-1} g'$ belongs to the local group of the smallest simplex of $Y$ that contains $x$. In particular, if $G = \widehat{G(\mathcal{Q})}$, then the space $D(Y,\varphi)$ is called the \textbf{universal cover} of $G(\mathcal{Q})$ with fundamental domain $Y$, and it is connected and simply connected (\cite{bridson2013metric}, Theorem II.12.20).
\end{defi}

\begin{rem} Let $G(\mathcal{Q})$ be a complex of groups and let $\varphi : G(\mathcal{Q}) \rightarrow G$ be a simple morphism to some group $G$, that is injective on the local groups. Then $G(\mathcal{Q})$ is developable. Indeed, let $Y$ be a geometric realisation of $\mathcal{Q}$ (i.e. the simplices of $Y$ form a poset that is isomorphic to $\mathcal{Q}$), then the space $X \coloneqq D(Y,\varphi)$ satisfies all the required properties of Definition \ref{DefiDev}.

In that case, the situation is quite simple: $X$ is a simplicial complex on which $G$ acts in such a way that the stabiliser of a simplex of the form $(1, \sigma)$ is precisely $G_{\sigma}$. Note that the stabiliser of a simplex of the form $(g, \sigma)$ is then exactly $g \cdot G_{\sigma} \cdot g^{-1}$.
\end{rem}

\section{The modified Deligne complex}

This section is dedicated to the modified Deligne complex. We will first recall its definition and then we will talk about its links and their geometry. The following definition is specific to Artin groups of dimension $2$. The definition in the more general case can be found in \cite{charney1995k}.

\begin{defi} \label{DefDeligne} Let $A_{\Gamma}$ be a $2$-dimensional Artin group of rank at least $3$, and let $\Delta_{\emptyset}$ be a $(|V(\Gamma)| - 1)$-simplex whose codimension $1$ faces are labeled as the set $\{\Delta_a \ | \ a \in V(\Gamma) \}$. In the barycentric subdivision of $\Delta_{\emptyset}$, we consider the following vertices:
\\$\bullet$ \underline{Type $0$:} $v_{\emptyset}$ is the vertex corresponding to $\Delta_{\emptyset}$.
\\$\bullet$ \underline{Type $1$:} for every $a \in S$, $v_a$ is the vertex corresponding to $\Delta_a$.
\\$\bullet$ \underline{Type $2$:} for every distinct $a, b \in S$, $v_{ab}$ is the vertex corresponding to $\Delta_a \cap \Delta_b$.
\\Let $K_{\Gamma}$ be the full subcomplex of the barycentric subdivision of $\Delta_{\emptyset}$ spanned by the vertex of type $0$, the vertices of type $1$, and the vertices of type $2$ that satisfy $m_{ab} < \infty$. Notice that $K_{\Gamma}$ is a geometric realisation of the poset $\mathcal{Q}$ of the vertices of type $0$, type $1$ and type $2$ for which $m_{ab} < \infty$. We define a complex of groups $G(\mathcal{Q})$ over $\mathcal{Q}$ in the following way. The local groups associated with $v_{\emptyset}$, $v_a$ and $v_{ab}$ are respectively $\{1\}$, $\langle a \rangle$ and $A_{ab}$. The natural inclusions of simplices
$$\Delta_a \cap \Delta_b \subseteq \Delta_a \subseteq \Delta_{\emptyset}$$
induce natural inclusions of the local groups
$$A_{ab} \supseteq \langle a \rangle \supseteq \{1\}$$
that define the maps $\psi_{\tau \sigma}$ of $G(\mathcal{Q})$. The simple morphism is the map $\varphi$ that associates to any vertex the natural inclusion of its local group into $A_{\Gamma}$. It follows from the definitions that $A_{\Gamma}$ is the fundamental group of $G(\mathcal{Q})$, whose universal cover is a $2$-dimensional space called the \textbf{modified Deligne complex} associated to $A_{\Gamma}$, and is denoted $D_{\Gamma}$.

Let us denote the edges of $K_{\Gamma}$ by $e_a$ if it connects $v_{\emptyset}$ and $v_a$, $e_{ab}$ if it connects $v_{\emptyset}$ and $v_{ab}$ and $e_{a,ab}$ if it connects $v_a$ and $v_{ab}$. Let us set $T_{ab}$ to be the unique triangle with vertices $v_{\emptyset}$, $v_a$ and $v_{ab}$. Note that $T_{ab}$ and $T_{ba}$ belong to $K_{\Gamma}$ whenever $m_{ab} < \infty$.

We will now try to give a better understanding of $K_{\Gamma}$ and $D_{\Gamma}$, notably in terms of angles and metrics. Consider a triangle $T_{ab} \subseteq K_{\Gamma}$. We recall the Moussong metric on $D_{\Gamma}$ as studied by Charney and Davis (\cite{charney1995k}). We first define the angles on $T_{ab}$ by the following:
\begin{align*}
\angle_{v_{ab}}(v_{\emptyset},v_a) &\coloneqq \frac{\pi}{2 \cdot m_{ab}} \\
\angle_{v_a}(v_{\emptyset},v_{ab}) &\coloneqq \frac{\pi}{2} \\
\angle_{_{\emptyset}}(v_a,v_{ab}) &\coloneqq \frac{\pi}{2} - \frac{\pi}{2 \cdot m_{ab}} 
\end{align*}
Notice that the angles of $T_{ab}$ add up to $\pi$, so that $T_{ab}$ is an euclidean triangle. We fix the length of every edge of the form $e_s$ to be $1$. The length of each edge of the form $e_{st}$ or $e_{s,st}$ can then be found by the law of sines. The metric on $K_{\Gamma}$ is the piecewise euclidean metric obtained by gluing the euclidean metrics from every triangle $T_{st}$. 

\begin{figure}[H]
\centering
\includegraphics[scale=0.5]{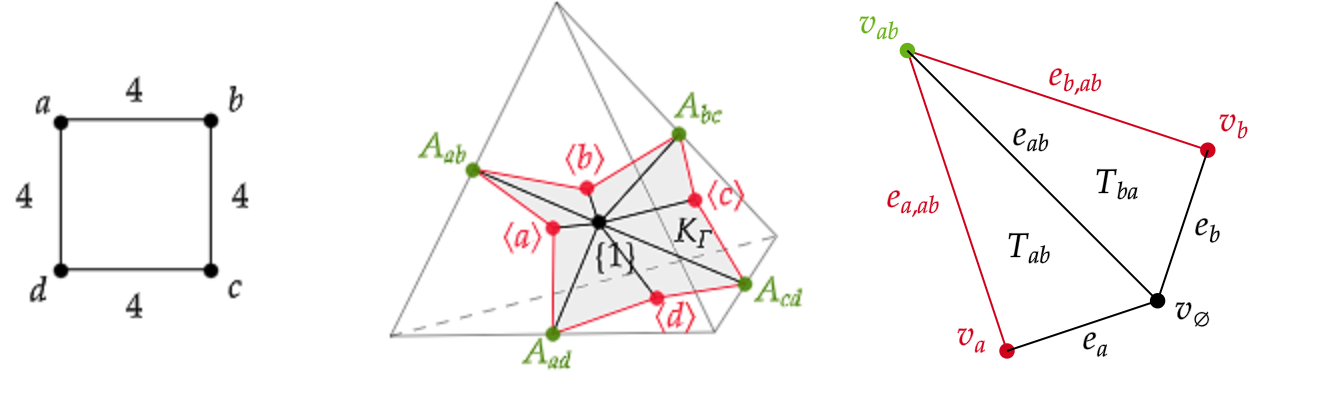}
\caption{\underline{On the left:} A graph $\Gamma$ defining an Artin group $A_{\Gamma}$.
\underline{In the centre:} The complex of groups $K_{\Gamma}$ with fundamental group $A_{\Gamma}$, and its local groups. It is obtained by gluing $8$ triangles, seen as a subset of a $3$-simplex.
\underline{On the right:} Example of notations for the part of $K_{\Gamma}$ corresponding to the triangles $T_{ab}$ and $T_{ba}$.
\\Vertices and edges have colors corresponding to the type of their local group (or stabiliser): black for the trivial group, red for an infinite cyclic group, and green for a dihedral Artin group.}
\label{FigModDelCpx}
\end{figure}

In light of (\cite{bridson2013metric}, Theorem II.12.18), the \textbf{modified Deligne complex} associated to $A_{\Gamma}$ can also be described as the space
$$D_{\Gamma} \coloneqq D(K_{\Gamma},\varphi) = \quotient{A_{\Gamma} \times K_{\Gamma}}{\sim},$$
where $(g,x) \sim (g',x') \Longleftrightarrow x = x'$ and $g^{-1} g'$ belongs to the local group of the smallest simplex of $K_{\Gamma}$ that contains $x$. Note that $A_{\Gamma}$ acts naturally on itself via left multiplication, and this action induces an action of $A_{\Gamma}$ on $D_{\Gamma}$ by simplicial morphisms with strict fundamental domain $K_{\Gamma}$.
\end{defi}

\begin{figure}[H]
\centering
\includegraphics[scale=1]{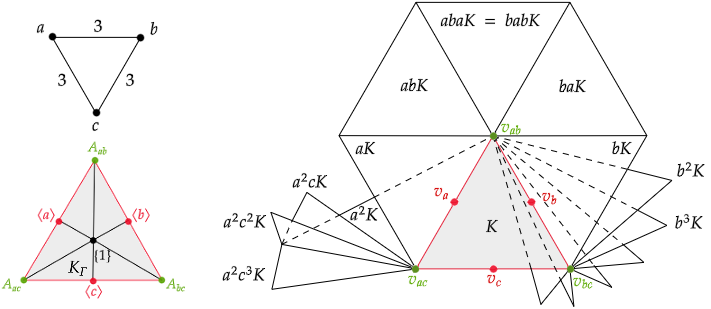}
\caption{\underline{On the top-left:} A graph $\Gamma$ defining an Artin group $A_{\Gamma}$.
\underline{On the bottom-left:} The complex of groups $K_{\Gamma}$ with fundamental group $A_{\Gamma}$, and its local groups.
\underline{On the right:} Part of the modified Deligne complex $D_{\Gamma}$. For drawing purposes, we wrote $K$ instead of $K_{\Gamma}$ and didn't draw the barycentric subdivision.
\\Vertices and edges have colors corresponding to the type of their local group (or stabiliser): black for the trivial group, red for an infinite cyclic group, and green for a dihedral Artin group.}
\end{figure}

\begin{rem} \label{RemDelCpx} (1) Let $\Gamma' \subseteq K_{\Gamma}$ be the graph defined by the vertices of type $1$ and $2$ and the edges $e_{a,ab}$ between them. Notice that  $K_{\Gamma}$ is just a cone over $\Gamma'$ with apex $v$, and that $\Gamma'$ is simply the barycentric subdivision of $\Gamma$. Thus, $K_{\Gamma}$ is the cone over the barycentric subdivision of $\Gamma$, with appropriate metric.
\\(2) Obtaining $A_{\Gamma}$ as the fundamental group of a complex of groups can be done is many different ways. The associated structure of complex of groups of $K_{\Gamma}$ described in Definition \ref{DefDeligne} is one of them, and we will call it the \textbf{usual} complex of groups associated with $A_{\Gamma}$. As we will see later, $A_{\Gamma}$ can be obtained as the fundamental group of larger complexes of groups (see Definition \ref{DefAugDeligne} for instance).
\end{rem}

\begin{defi} Let $X$ be a simplicial complex and let $v \in X$. The \textbf{link} of $v$ in $X$ is the set $Lk_X(v)$ of simplices of $X$ that are disjoint from $v$ but that belong to simplices of $X$ that contain $v$. We endow $Lk_X(v)$ with the following metric:
$$\forall \text{ simplices } \sigma, \tau \in Lk_X(v), \forall \ x \in \sigma, \forall \ y \in \tau, \ \ d_{Lk_X(v)}(x,y) \coloneqq \angle_v(x,y),$$
where $\angle_v(x,y)$ is the angular distance between $x$ and $y$ (in the sense of \cite{bridson2013metric}).
\end{defi}

\noindent \begin{rem} \label{RemLinks} Following the previous definitions, a natural question to ask is what do the links of vertices of $D_{\Gamma}$ look like ? In light of (\cite{bridson2013metric}, Construction II.12.24), the link $Lk_{D_{\Gamma}}(v)$ around a vertex $v \in K_{\Gamma}$ only depends on the development of the local groups around $v$. More specifically, the link $Lk_{D_{\Gamma}}(v)$ is isomorphic to the development $D(Lk_{K_{\Gamma}}(v),(\psi_v)_{e \centerdot})$ of the link $Lk_{K_{\Gamma}}(v)$ along the natural inclusion maps $(\psi_v)_{e \centerdot} : G_e \hookrightarrow G_v$, where $e$ is an edge from $v$ to $Lk_{D_{\Gamma}}(v)$ and $e \centerdot \coloneqq e \cap Lk_{K_{\Gamma}}(v)$. In particular, we can give a more precise geometric description of the links of vertices in $D_{\Gamma}$:
\\$\bullet$ \underline{Type $0$:} $Lk_{D_{\Gamma}}(v_{\emptyset})$ is the development of $Lk_{K_{\Gamma}}(v_{\emptyset})$ over the trivial maps $(\psi_{v_{\emptyset}})_{e_a}: \{1\} \hookrightarrow \{1\}$ and $(\psi_{v_{\emptyset}})_{e_{ab}}: \{1\} \hookrightarrow \{1\}$. Notice that $Lk_{K_{\Gamma}}(v)$ is the graph $\Gamma'$ from Remark \ref{RemDelCpx}, and hence $Lk_{D_{\Gamma}}(v_{\emptyset})$ is just the barycentric subdivision of $\Gamma$. By construction, the lengths of edges in $Lk_{D_{\Gamma}}(v_{\emptyset})$ are given by
$$\ell(e_{a,ab}) =  \angle_{v_{\emptyset}}(v_a,v_{ab}) = \frac{\pi}{2} - \frac{\pi}{2 \cdot m_{ab}}.$$
$\bullet$ \underline{Type $1$:} $Lk_{D_{\Gamma}}(v_a)$ is the development of $Lk_{K_{\Gamma}}(v_a)$ over the maps $(\psi_{v_a})_{e_a}: \{1\} \hookrightarrow \langle a \rangle$ and $(\psi_{v_a})_{e_{a,ab}}: \langle a \rangle \hookrightarrow \langle a \rangle$. It is not hard to see that $Lk_{K_{\Gamma}}(v_a)$ is just a $n_a$-pod centered at $v_{\emptyset}$, where $n_a \coloneqq | \{b \in V(\Gamma) \backslash \{a\} \ | \ m_{ab} < \infty \} |$. In particular, $Lk_{D_{\Gamma}}(v_a)$ is the quotient $\quotient{Lk_{K_{\Gamma}}(v_a) \times \langle a \rangle}{\sim},$ where $(x,a^n) \sim (y,a^m)$ if and only if either $x = y = v_{ab}$ for some $b \in V(\Gamma) \backslash \{a\}$ with $m_{ab} < \infty$ or $x = y$ and $n = m$. Notice that by construction, every edge $e_{ab}$ has length $\angle_{v_a}(v_{\emptyset},v_{ab}) = \pi/2$ in $Lk_{D_{\Gamma}}(v_a)$.
\\$\bullet$ \underline{Type $2$:} $Lk_{D_{\Gamma}}(v_{ab})$ is the development of $Lk_{K_{\Gamma}}(v_{ab})$ over the three maps $(\psi_{v_{ab}})_{e_{ab}}: \{1\} \hookrightarrow A_{ab}$, $(\psi_{v_{ab}})_{e_{a,ab}}: \langle a \rangle \hookrightarrow A_{ab}$ and $(\psi_{v_{ab}})_{e_{b,ab}}: \langle b \rangle \hookrightarrow A_{ab}$. The link $Lk_{K_{\Gamma}}(v_{ab})$ is simply a tree $T_0$ that consists of the two edges $e_a$ and $e_b$. Consider the Bass-Serre tree $T$ over $T_0$ and its associated local groups. The development of $T_0$ over the previously described maps is just the quotient of $T$ by $\langle \langle \underbrace{aba\cdots}_{m_{ab}} = \underbrace{bab\cdots}_{m_{ab}} \rangle \rangle$, because the previous maps inject into
$$A_{ab} \cong \quotient{F_{ab}}{\langle \langle \underbrace{aba\cdots}_{m_{ab}} = \underbrace{bab\cdots}_{m_{ab}} \rangle \rangle}.$$
Notice by construction that the lengths of edges in $Lk_{D_{\Gamma}}(v_{ab})$ are given by
$$\forall s \in \{a,b\}, \ \ell(e_s) =  \angle_{v_{ab}}(v_{\emptyset},v_s) = \frac{\pi}{2 \cdot m_{ab}}.$$

\begin{figure}[H]
\centering
\includegraphics[scale=1.1]{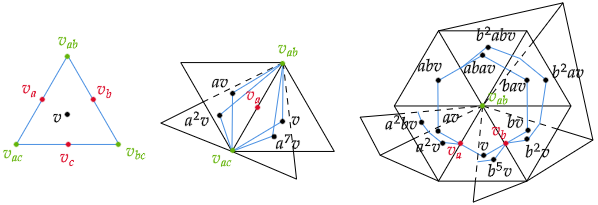}
\caption{Part of the links of the vertices of type 0, 1 and 2 respectively, from left to right. The links are drawn in blue.
For drawing purposes, we wrote $v$ instead of $v_{\emptyset}$.}
\end{figure}
\end{rem}

\noindent Using the description of links seen in Remark \ref{RemLinks}, Charney and Davis proved the following:

\begin{thm} \label{XCat0} \textbf{(\cite{charney1995k}, Proposition 4.4.5)} Let $A_{\Gamma}$ be a $2$-dimensional Artin group of rank at least $3$. Then its modified Deligne complex $D_{\Gamma}$ is CAT(0).
\end{thm}

\section{Links of vertices of type 2}

The goal of this section is to get a better understanding of the links $Lk_{D_{\Gamma}}(v_{ab})$ of vertices of type $2$ in $D_{\Gamma}$, and ultimately to prove Proposition E that we recall here:

\begin{prop*} Let $A_{\Gamma}$ be a $2$-dimensional Artin group of rank at least $3$ with modified Deligne complex $D_{\Gamma}$. Suppose that there exists a vertex $v_{ab} \in D_{\Gamma}$ whose stabiliser $A_{ab}$ has coefficient $3 \leq m_{ab} < \infty$. Then:
\\(1) The orbits of $A_{ab}$ on $Lk_{D_{\Gamma}}(v_{ab})$ are unbounded.
\\(2) More precisely, the orbits of $\langle g \rangle$ on $Lk_{D_{\Gamma}}(v_{ab})$ are quasi-isometrically embedded if and only if $g \in A_{ab}$ is not trivial, nor the conjugate of a power of a standard generator $a$ or $b$.
\end{prop*}

In particular, we would like to prove that if $A_{\Gamma}$ has a coefficient $3 \leq m_{ab} < \infty$, then the first condition of Theorem D is satisfied for $v_{ab}$ (Proposition E.(1)). We also prove results of independent interests, classifying which elements $g \in A_{ab}$ generate subgroups $\langle g \rangle$ that act on $Lk_{D_{\Gamma}}(v_{ab})$ with quasi-isometrically embedded orbits (Proposition E.(2)). This result was used by Hagen-Martin-Sisto in proving that extra-large type Artin groups are virtually hierarchically hyperbolic (\cite{hagen2021extra}).

\subsection{Reformulating Proposition E in terms of syllabic lengths}

The goal of this section is to reformulate Proposition E into a more accessible problem (see Proposition \ref{PropEquivToPropE}). We begin with the following definition, that will be useful throughout all the section: 

\begin{defi} Let $G$ be a group with generating set $\mathcal{A}$, and let $\varphi : F_{\mathcal{A}} \twoheadrightarrow G$ be the natural surjection from the free group over $\mathcal{A}$ onto $G$.
\\$\bullet$ Every word $w \in F_{\mathcal{A}}$ can be written uniquely as $w = a_1^{r_1} \cdots a_n^{r_n}$, assuming $a_i \in \mathcal{A}$, $a_i \neq a_{i+1}$ and $r_i \in \mathbf{Z} \backslash \{0\}$. Then the \textbf{syllabic length} of $w$ is $\ell_{\mathcal{S}}(w) \coloneqq n$.
\\$\bullet$ For every element $g \in G$ we define the \textbf{syllabic length} of $g$ as $\ell_{\mathcal{S}}(g) \coloneqq \min \{ \ell_{\mathcal{S}}(w) \ | \ \varphi(w) = g \}$.
\end{defi}

Recall that in the modified Deligne complex $D_{\Gamma}$ associated with an Artin group $A_{\Gamma}$, the stabilisers of vertices of type 2 (ex: $v_{ab}$) are dihedral Artin groups (ex: $A_{ab}$). The following lemma makes a connection between the syllabic length of elements $g \in A_{ab}$ and the distances in the link $Lk_{D_{\Gamma}}(v_{ab})$.

\begin{lemma} \label{LemmaSyl} Let $\gamma$ be a path in $Lk_{D_{\Gamma}}(v_{ab})$ joining $v_{\emptyset}$ and $g v_{\emptyset}$ for some $g \in A_{ab}$, and suppose that the edges of $\gamma$ are $e_a, \ a^{n_1} e_a, \ a^{n_1} e_b, \ a^{n_1} b^{n_2} e_b, \cdots, \ a^{n_1} b^{n_2} \cdots x^{n_k} e_x$, in that order, where $x \in \{a,b\}$ and $n_i \in \mathbf{Z} \backslash\{0\}$. Let now $w \coloneqq a^{n_1} b^{n_2} \cdots x^{n_k}$. Then $\ell(\gamma) = \frac{\pi}{m_{ab}} \cdot \ell_{\mathcal{S}}(w)$. Furthermore, $d_{Lk_{D_{\Gamma}}(v_{ab})}(v_{\emptyset}, g v_{\emptyset}) = \frac{\pi}{m_{ab}} \cdot \ell_{\mathcal{S}}(g)$.
\end{lemma}

\noindent \textbf{Proof:} First of all, recall that the local group at $v_a$ is $\langle a \rangle$, and hence the set of edges of $Lk_{D_{\Gamma}}(v_{ab})$ meeting at $v_a$ is $\{a^k e_a \ | \ k \in \mathbf{Z} \}$. This proves that $e_a$ and $a^{n_1} e_a$ are indeed consecutive to one another, meeting at $v_a$. Of course, $a^{n_1} e_a$ and $a^{n_1} e_b$ are also consecutive to one another, meeting at $a^{n_1} v_{\emptyset}$. A similar argument shows that the edges in the statement of the lemma consecutively meet each others. As $A_{\Gamma}$ acts by isometries on $Lk_{D_{\Gamma}}(v_{ab})$, it is clear that the length of every edge of $\gamma$ is either $\ell(e_a)$ or $\ell(e_b)$, both of which turn out to be equal to $\frac{\pi}{2 \cdot m_{ab}}$. Because the number of edges in $\gamma$ is precisely $2 \cdot \ell_{\mathcal{S}}(w)$, we get $\ell(\gamma) = \frac{\pi}{m_{ab}} \cdot \ell_{\mathcal{S}}(w)$

Notice that every path $\gamma$ joining $v_{\emptyset}$ and $g v_{\emptyset}$ corresponds to a word $w$ that satisfies $\varphi(w) = g$, where $\varphi: F_{ab} \twoheadrightarrow A_{ab}$ is the natural projection. The distance between $v_{\emptyset}$ and $g v_{\emptyset}$ is the length of the shortest of these paths, hence
$$d_{Lk_{D_{\Gamma}}(v_{ab})}(v_{\emptyset}, g v_{\emptyset}) = \min \{ \frac{\pi}{m_{ab}} \cdot \ell_{\mathcal{S}}(w) \ | \ \varphi(w) = g \} = \frac{\pi}{m_{ab}} \cdot \ell_{\mathcal{S}}(g).$$
\hfill\(\Box\)

One important consequence of the previous lemma is that we can reformulate Proposition E in terms of syllabic lengths of elements in the local group of a vertex of type $2$. We will prove Proposition E by proving the following equivalent proposition:

\begin{prop} \label{PropEquivToPropE} Let $A_{\Gamma}$ be a $2$-dimensional Artin group of rank at least $3$ with modified Deligne complex $D_{\Gamma}$. Suppose that there exists a vertex $v_{ab} \in D_{\Gamma}$ whose stabiliser $A_{ab}$ has coefficient $3 \leq m_{ab} < \infty$. Then:
\\(1) The set $\{ \ell_{\mathcal{S}}(g) \ | \ g \in A_{ab} \}$ is unbounded.
\\(2) More precisely, the syllabic length $\ell_{\mathcal{S}}(g^n)$ grows linearly in $n$ if and only if $g \in A_{ab}$ is not trivial, nor the conjugate of a power of a standard generator $a$ or $b$.
\end{prop}

\begin{rem} Recall that we say that a sequence $\{u_n \}_{n \geq 0}$ grows linearly in $n$ if they are constants $B \geq A > 0$ and $C \geq 0$ such that for any $n \geq 0$ we have
$$A n - C \leq u_n \leq B n + C.$$
\end{rem}

The next lemma shows that Proposition \ref{PropEquivToPropE}.(1), and thus Proposition E.(1), are satisfied. This result will be very useful in the proof of Theorem A. It shows that if $A_{\Gamma}$ has a coefficient $m_{ab} \geq 3$, then the vertex $v_{ab}$ satisfies the first hypothesis of Theorem D.

\begin{lemma} \label{LemmaInfWords} Consider an Artin group $A_{ab}$ with coefficient $3 \leq m_{ab} \leq \infty$. Then
$$\{\ell_{\mathcal{S}}(g) \ | \ g \in A_{ab} \} \text{ is unbounded}.$$
\end{lemma}

\noindent \textbf{Proof:} It is known that the quotient $\bar{A}_{ab}$ of $A_{ab}$ by its centre is virtually isomorphic to the free group $F_m$, for some $m \geq 2$ (\cite{brady2000three},\cite{ciobanu2020equations}). In particular, $\bar{A}_{ab}$ is acylindrically hyperbolic. Suppose now that there exists a constant $N \geq 0$ such that for every $g \in A_{ab}$, one has $\ell_{\mathcal{S}}(g) < N$, and assume without loss of generality that $N$ is even. This means that $A_{ab} = \langle a \rangle \langle b \rangle \cdots \langle a \rangle \langle b \rangle$ (where the product has $N$ terms). In particular, $\bar{A}_{ab} = \bar{\langle a \rangle} \bar{\langle b \rangle} \cdots \bar{\langle a \rangle} \bar{\langle b \rangle}$. We can now use (\cite{osin2016acylindrically}, Proposition 1.7) to see that one of $\bar{\langle a \rangle}$ or $\bar{\langle b \rangle}$ must be acylindrically hyperbolic, which is impossible, as they are cyclic subgroups of $\bar{A}_{ab}$. Therefore, $\{\ell_{\mathcal{S}}(g), g \in A_{ab} \}$ is unbounded.
\hfill\(\Box\)
\bigskip

\noindent \textbf{Strategy:} The goal of the rest of this section is to understand more those links of the form $Lk_{D_{\Gamma}}(v_{ab})$, i.e. the links of vertices of type $2$ in $D_{\Gamma}$. In particular, we will be able through a more precise analysis of these links to prove Proposition \ref{PropEquivToPropE}.(2), and thus Proposition E.(2).

We now set for the rest of this section $A_{ab}$ to be a dihedral Artin group with coefficient $3 \leq m < \infty$. One can easily see that for any element $g \in A_{ab}$, the syllabic length $\ell_{\mathcal{S}}(g^n)$ is always bounded above by a linear function, such as $\ell(w) \cdot n$ for instance, where $w$ is any word representing $g$ and $\ell(\cdot)$ is the usual length function on words. Therefore we will only focus on finding a linear lower bound for $\ell_{\mathcal{S}}(g^n)$.

Our approach is mostly geometric: we study the action of $A_{ab}$ on a graph $\widehat{T}$ (see Definition \ref{DefiTHat}). In particular, we show that the distance of translation induced by an element $g \in A_{ab}$ gives a lower bound on the syllabic length of $g$ (see Lemma \ref{LemmaDistanceAndSyllabicLengthInHatT}). It then follows immediately that any element $g \in A_{ab}$ that acts loxodromically on $\widehat{T}$ is such that $\ell_{\mathcal{S}}(g^n)$ admits a linear lower bound in $n$, giving Proposition \ref{PropEquivToPropE}.(2) for such elements. It then feels natural to want to determine which elements act loxodromically on $\widehat{T}$. This will be achieved in Lemma \ref{LemmaEllipticOfHatT}.

It remains to study the elements that do not act loxodromically on $\widehat{T}$. They all act elliptically and come in two forms: the elements that are conjugate to powers of a standard generator (modulo an element of the centre), and the (non-trivial) elements which admit powers that belong to the centre of $A_{ab}$. When their "central part" is trivial, the elements $g$ of the first kind are easily shown to satisfy $\ell_{\mathcal{S}}(g^n) \leq K_g$ for a constant $K_g$ that does not depend on $n$. However, the elements of the first kind that don't have a trivial central part and the elements of the second kind have a different behaviour. As will be recalled later, the centre of $A_{ab}$ only contains powers of the Garside element of $A_{ab}$, which motivates a more in-depth study of the syllabic length of such powers. The method that we use for that last point is more algebraic, and rely on a more direct study of the syllabic length of words, notably using the Garside normal form of elements. As a consequence, we will show that the remaining elliptic elements $g$ are such that the syllabic length $\ell_{\mathcal{S}}(g^n)$ also admits a linear lower bound in $n$ (see Lemma \ref{LemmaLowerBoundForDelta}). Alltogether, this will conclude the proof of Proposition \ref{PropEquivToPropE}.

\subsection{The action on $\widehat{T}$.}

Our first goal is to define a tree $T$ on which $A_{ab}$ acts nicely with a trivial action of the centre $Z(A_{ab})$. This will be done throughout the next definitions and lemmas. Let us first introduce few notations and recall notions about normal forms and centres in dihedral Artin groups.
\bigskip

\noindent \textbf{Notations:} $\bullet$ For two words $u, u' \in F_{ab}$ representing the same element $g \in A_{ab}$, we will simply write $u \equiv u'$ instead of $\varphi(u) = \varphi(u')$. Similarly, for $g \in A_{ab}$, we will write $u \equiv g$ instead of $\varphi(u) = g$.
\\$\bullet$ We will write $(a,b;k)$ to denote the alternating sequence of the letters $a$ and $b$, starting with $a$ and of length $k$, and we will write $\Delta_a$ and $\Delta_b$ to describe the words $(a,b;m) \in F_{ab}$ and $(b,a;m) \in F_{ab}$ respectively. More explicitely,
$$\Delta_a \coloneqq \underbrace{aba \cdots}_{m \text{ terms}} \ \text{ and } \ \Delta_b \coloneqq \underbrace{bab \cdots}_{m \text{ terms}}.$$
\\$\bullet$ For a word $u \in F_{ab}$, we denote by $\bar{u}$ the word obtained from $u$ by replacing every $a^n$ by $b^n$ and every $b^n$ by $a^n$. Moreover, we will denote by $\widetilde{u}$ the element
$$\widetilde{u} \coloneqq \left\{\begin{array}{ll}
u & \text{if } m \text{ is even} \\
\bar{u} & \text{if } m \text{ is odd}
\end{array} \right.$$
One can easily notice that for any word $u \in F_{ab}$, we have
$\Delta^{\pm1} \cdot u \equiv \widetilde{u} \cdot \Delta^{\pm1}$.

\begin{defi} \label{DefiGarside} For a dihedral Artin group $A_{ab}$ with coefficient $3 \leq m_{ab} < \infty$, the \textbf{Garside element} is the element $\Delta \in A_{ab}$ defined by
$$\Delta \equiv \Delta_a \equiv \Delta_b. \ \ (*)$$
A strict non-trivial subword of $\Delta_a$ or of $\Delta_b$ is called an \textbf{atom}. It is a standard result (\cite{garside1969braid}, \cite{elrifai1994algorithms}, \cite{dehornoy1999gaussian}) that for every element $g \in A_{ab}$, there is a word $\operatorname{Gars}(g) \in F_{ab}$ called the \textbf{Garside normal form} of $g$ that satisfies $\operatorname{Gars}(g) \equiv g$ and such that one can write
$$\operatorname{Gars}(g) = u_1 \cdots u_n \cdot W,$$
where the $u_i$ are atoms such that the last letter of each $u_i$ matches with the first letter of $u_{i+1}$, and where $W \equiv \Delta^N$ for some $N \in \textbf{Z}$ is a product of terms of the form $\Delta_a^{\pm1}$ and $\Delta_b^{\pm1}$. This word is not unique, however the atoms of the above decomposition are uniquely defined, and so is $N$.
\end{defi}

\noindent At last, we recall that the centre $Z(A_{ab})$ of $A_{ab}$ was described in \cite{brieskorn1972artin}, and takes the form
$$Z(A_{ab}) = \begin{cases} \langle \Delta \rangle \text{ if } m \text{ is even}, \\
\langle \Delta^2 \rangle \text{ if } m \text{ is odd}. \end{cases}$$

\noindent We now come back to constructing the desired space. The space that we first define is due to (\cite{bestvina1999non}, Section 2.1). One can also recover an equivalent definition by quotient of the space described in (\cite{mccammond2010combinatorial}, Figure 6).

\begin{figure}[H]
\centering
\includegraphics[scale=0.5]{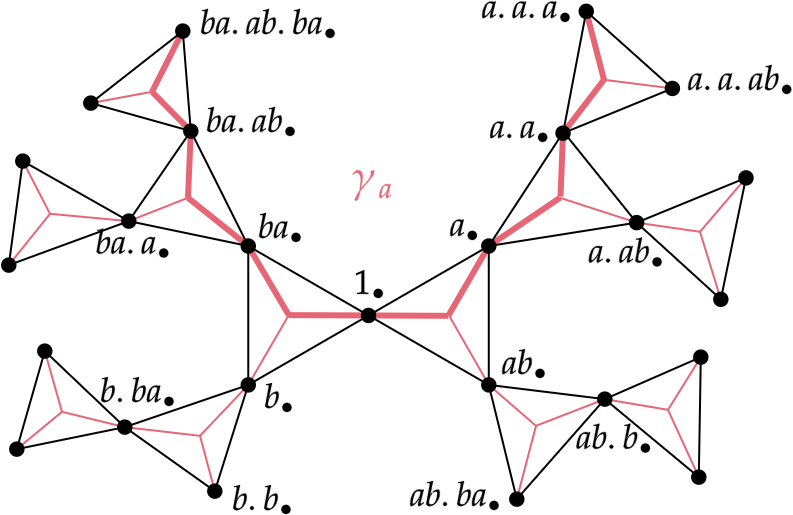}
\caption{Let $m \coloneqq 3$. \underline{In black:} Part of the graph $Y$ with its set of vertices $V$. \underline{In pink:} Part of the tree $T$, that is a deformation retract of $Y$. The axis $\gamma_{a} \subseteq T$ is drawn with the thicker line.}
\label{FigureTree}
\end{figure}

\begin{defi} \label{DefiY} Let us consider the graph $Y$ defined by the following (see Figure \ref{FigureTree}):
\bigskip

\noindent \underline{Vertices:} The vertex set of $Y$ is the set of cosets
$$V \coloneqq \quotient{A_{ab}}{\langle \Delta \rangle} = \{ g \langle \Delta \rangle \ | \ g \in A_{ab} \}.$$
A convenient representative for a vertex $g \langle \Delta \rangle$ is the product of the atoms of the Garside normal form of $g$. This representative is the unique that is in Garside normal form yet does not contain any subword of the form $\Delta_x^{\pm 1}$ for some $x \in \{a, b\}$. We will denote it $g_{\bullet}$. In this setup, we can see $V$ as the set $\{ g_{\bullet} \ | \ g \in A_{ab} \}$.
\bigskip

\noindent \underline{Simplices:} For every collection $g_{1 \bullet}, \cdots, g_{k \bullet}$ of vertices, the set $\{g_{1 \bullet}, \cdots, g_{k \bullet} \}$ spans a $k$-simplex if and only if for all $i, j \in \{1, \cdots, k\}$, there is an atom $x$ such that $g_{i \bullet} \cdot x = g_{j \bullet}$ or $g_{j \bullet} \cdot x = g_{i \bullet}$. Note that because atoms are subwords of $\Delta_a$ or $\Delta_b$, every $k$-simplex is contained in a maximal $m$-simplex, where $m$ is the coefficient of $A_{ab}$ (see Figure \ref{FigureTree}).
\bigskip

\noindent The group $A_{ab}$ acts naturally on $V$: if $h \in A_{ab}$ and $g \langle \Delta \rangle \in V$, then $h \cdot g \langle \Delta \rangle \coloneqq hg \langle \Delta \rangle$. This action extends to a simplicial and cocompact action of $A_{ab}$ on $Y$, that is transitive on the vertices (see \cite{bestvina1999non}). Note that $Z(A_{ab})$ acts trivially on $V$, and thus on $Y$.
\end{defi}

\begin{defi} \label{DefiT} We define a new graph $T$ by the following. The set of vertices of $T$ is the union of two sets: the set $V$ of vertices of $Y$, and the set $V'$ of maximal simplices of $Y$ (i.e. the $m$-simplices). Then, we put an edge between a vertex $g_{\bullet} \in V$ and a vertex $\{g_{1 \bullet}, \cdots, g_{m \bullet} \} \in V'$ if and only if $g_{\bullet} \in \{g_{1 \bullet}, \cdots, g_{m \bullet} \}$). Note that $T$ can naturally be seen as a subspace of $Y$ (see Figure \ref{FigureTree}).
\end{defi}

\begin{lemma} \label{LemmaTTree} The graph $T$ is a tree, and the stabiliser of any edge $e \subseteq T$ is precisely $Z(A_{ab})$.
\end{lemma}

\noindent \textbf{Proof:} The atoms in the Garside normal form of any element are unique, and this gives $Y$ a structure of tree of $m$-simplices (see Figure \ref{FigureTree}), where the $Y$-distance between any vertex $g_{\bullet}$ and $1_{\bullet}$ is precisely the number of atoms in $g_{\bullet}$ (note that the $T$-distance is twice that amount). The reason $m$-simplices appear is because once given a non-trivial vertex $g_{\bullet}$, there are $(m-1)$ different ways one can add an atom on the right side of $g_{\bullet}$ (assuming this atom starts with a letter that differs from the last letter of $g_{\bullet}$. In particular, $Y$ retracts on a tree described in Figure \ref{FigureTree}, and that tree is precisely $T$.

Let now $e$ be any edge of $T$. Because the action is transitive on the vertex set $V$, we may as well assume that $e$ contains $1_{\bullet}$. The other vertex of $e$ corresponds to one of the two simplices $S_a \coloneqq \{ 1_{\bullet}, a_{\bullet}, \cdots, (a,b;m-1)_{\bullet} \}$ or $S_b \coloneqq \{ 1_{\bullet}, b_{\bullet}, \cdots, (b,a;m-1)_{\bullet} \}$. Let now $g \in A_{ab}$ and suppose that $g \cdot e = e$. Then in particular $g$ fixes $1_{\bullet}$, so we have $g \cdot \langle \Delta \rangle = \langle \Delta \rangle$, and thus $g \in \langle \Delta \rangle$. If $m$ is even, we are done. If $m$ is odd, it is enough to show that $\Delta$ does not fix $e$, which is clear because it sends $S_a$ onto $S_b$ and vice versa.
\hfill\(\Box\)

\begin{rem} \label{RemValence} The valence of a vertex $v$ of $T$ is easy to determine. If $v \in V$, then $v$ belongs to exactly two $m$-simplices of $Y$, so the valence of $v$ in $T$ is $2$. If $v \in V'$, then $v$ corresponds to a $m$-simplex of $Y$, hence is connected to exactly $m$ vertices of $Y$, and its valence in $T$ is $m$.
\end{rem}


\begin{lemma} \label{LemmaEllipticOfT} The elements of $A_{ab}$ acting elliptically on $T$ are precisely the elements $g \in A_{ab}$ for which there exists an $N \neq 0$ such that $g^N \in Z(A_{ab})$. All the other elements act loxodromically.
\end{lemma}

\noindent \textbf{Proof:} Suppose first that $g^N \in Z(A_{ab})$. Then $g^N$ acts trivially on $T$, so $g$ has finite orbits. In particular these orbits are bounded, so $g$ acts elliptically.

Suppose now that $g$ acts elliptically. Then $g$ fixes a vertex $v$ of $T$. By Remark \ref{RemValence}, $v$ has at most $m$ neighbours, so $g^{m !}$ fixes the neighbourhood of $v$. In particular, $g^{m!}$ fixes a vertex $h \langle \Delta \rangle$ of $Y$ (either $v$, or a vertex in its neighbourhood). The equation $g^{m !} \cdot h \langle \Delta \rangle = h \langle \Delta \rangle$ gives $g^{m !} = h \Delta^K h^{-1}$ for some $K \in \textbf{Z}$. We obtain $g^{2 m !} = h \Delta^{2K} h^{-1} = \Delta^{2K} \in Z(A_{ab})$.

All other elements act loxodromically because we have a simplicial isometric action on a tree.
\hfill\(\Box\)

Recall that we are interested in studying the syllabic length of elements of $A_{ab}$ relatively to the standard generators $a$ and $b$, and in reducing the problem of syllabic lengths to a problem of distances in our space. Unfortunately, one can travel an arbitrary large distance in $T$ using a single syllable, because the generators $a$ and $b$ act loxodromically on $T$ (see Figure \ref{FigureTree}). To deal with that problem, we decide to cone-off such axes, and their translates:

\begin{defi} \label{DefiTHat} Let $s \in \{a, b\}$. We denote by $\gamma_{s}$ the axis of $s$ in $T$, i.e. the bi-infinite geodesic line going through all the vertices of the form $s^k \langle \Delta \rangle$ for all $k \in \textbf{Z}$. Let us now define a graph $\widehat{T}$ as the cone-off of the tree $T$ along the family of axes $h \cdot \gamma_s$, for all $h \in A_{ab}$ and $s \in \{a, b\}$. More precisely:
\\$\bullet$ Start with $T$, and add a new vertex $v_{h, s}$ for every axis of the form $h \cdot \gamma_s$, for all $h \in A_{ab}$ and $s \in \{a, b\}$. We only add one vertex if two axes define the same line, even if they go in opposite directions.
\\$\bullet$ Connect every vertex $v_{h, s}$ to every vertex of the corresponding axis $h \cdot \gamma_s$.
\end{defi}

\noindent The following lemma justify the study of the cone-off $\widehat{T}$, as it gives a lower bound on the syllabic length of an element of $A_{ab}$ in terms of distances in $\widehat{T}$.

\begin{lemma} \label{LemmaDistanceAndSyllabicLengthInHatT} Let $g \in A_{ab}$. Then
$$d_{\widehat{T}}(1_{\bullet}, g \cdot 1_{\bullet}) \leq 2 \cdot \ell_{\mathcal{S}}(g).$$
\end{lemma}

\noindent \textbf{Proof:} The argument is similar to that of Lemma \ref{LemmaSyl}. Let $k \coloneqq \ell_{\mathcal{S}}(g)$. By the triangle inequality, it is enough to prove that there is a sequence of vertices
$$v_0 \coloneqq 1_{\bullet}, \ v_1, \ \cdots, \ v_{k-1}, \ v_k \coloneqq g \cdot 1_{\bullet}$$
such that $d_{\widehat{T}}(v_i, v_{i+1}) \leq 2$. To do so, let $a^{n_1} b^{n_2} \cdots x^{n_k}$ be a word representing $g$, where $x \in \{a, b\}$ and $n_i \in \textbf{Z} \backslash \{0\}$, and let
$$g_i \coloneqq a^{n_1} b^{n_2} \cdots y_i^{n_i},$$
where $y_i \in \{a, b\}$ is the appropriate letter. Consider now the vertices $v_i$ defined by $v_i \coloneqq g_i \cdot 1_{\bullet}$, so that $v_0 = 1_{\bullet}$ and $v_k = g \cdot 1_{\bullet}$. Because $A_{ab}$ acts on $\widehat{T}$ by isometries, we have for any $0 \leq i < k$
$$d_{\widehat{T}}(v_i, v_{i+1}) = d_{\widehat{T}}(g_i \cdot 1_{\bullet}, g_{i+1} \cdot 1_{\bullet}) = d_{\widehat{T}}(1_{\bullet}, y_{i+1}^{n+1} \cdot 1_{\bullet}).$$
Note that $y_{i+1}^{n+1}$ is just a power of a standard generator $y_{i+1} \in \{a, b\}$, which means $1_{\bullet}$ and $y_{i+1}^{n+1} \cdot 1_{\bullet}$ both belong to the axis $\gamma_{y_{i+1}}$. By definition of $\widehat{T}$, such vertices lie within distance $2$ of each others. It follows that $d_{\widehat{T}}(v_i, v_{i+1}) \leq 2$.
\hfill\(\Box\)
\bigskip

As explained in the strategy of this section, the previous lemma immediately gives Proposition \ref{PropEquivToPropE}.(2) for elements of $A_{ab}$ acting loxodromically on $\widehat{T}$. The goal of the next lemma is to classify these elements:

\begin{lemma} \label{LemmaEllipticOfHatT} The elements of $A_{ab}$ acting elliptically on $\widehat{T}$ are precisely the elements $g \in A_{ab}$ satisfy one of the following:
\\(1) $g = h \cdot s^N \cdot h^{-1} \cdot W$ for some $N \in \textbf{Z}$, $s \in \{a, b\}$ and $W \in Z(A_{ab})$;
\\(2) There exists an $N \neq 0$ such that $g^N \in Z(A_{ab})$.
\\All the other elements act loxodromically.
\end{lemma}

\noindent \textbf{Proof:} Let $g \in A_{ab}$. If $g$ satisfies $(2)$, then it already acts elliptically on $T$ by Lemma \ref{LemmaEllipticOfT}, so it acts elliptically on $\widehat{T}$ too. If $g$ satisfies $(1)$, it is not hard to see that $g$ fixes the vertex $h \cdot \gamma_{s}$, hence acts elliptically on $\widehat{T}$. We now suppose that $g$ does not satisfy any of these properties. We already know by Lemma \ref{LemmaEllipticOfT} that $g$ acts loxodromically on $T$, with an axis that we call $\gamma_g$. We begin by stating the following "small cancellation" claim, which gives the desired result for $g$. Then, we proceed on proving that $g$ actually satisfies the hypotheses of the claim.
\bigskip

\noindent \textbf{Claim:} Suppose that there exists a $K > 0$ such that for every $h \in A_{ab}$ and every $s \in \{a, b\}$, the subtree $\gamma$ defined by $\gamma \coloneqq \gamma_g \cap h \cdot \gamma_s$ has diameter at most $K$. Then $g$ acts loxodromically on $\widehat{T}$.
\bigskip

\noindent \textbf{Proof of the Claim:} Since $g$ acts loxodromically on $T$, it is enough to show that there is a constant $C > 0$ such that for all vertices $x, y \in \gamma_g$, we have
$$C \cdot d_T(x,y) \leq d_{\widehat{T}}(x,y). \ \ (*)$$
Let $\gamma_{x,y}^T$ be the (unique) geodesic connecting $x$ and $y$ in $T$, and let $M$ be the minimal number of axes of the form $h \cdot \gamma_s$ required to cover all edges of $\gamma_{x,y}^T$ completely. Let also $D \coloneqq d_{\widehat{T}}(x,y)$. Since every edge has length $1$, this means we can reach $x$ from $y$ by using $D$ edges $e_1, \cdots, e_D$ of $\widehat{T}$. Let $x_0, \cdots, x_D$ be the vertices these edges go through (in that order), and let $x_{r_0}, \cdots, x_{r_{D'}}$ be the subset of the above vertices corresponding to those belonging to $T$ (with $r_0 < \cdots < r_{D'}$). Then the vertices $x_{r_i}$ and $x_{r_{i+1}}$ always belong to a common axis. Indeed, if $r_{i+1} - r_i = 1$, then the two vertices are the two endpoints of a common edge of $T$. On the other hand, if $r_{i+1} - r_i \geq 2$, then there is at least one vertex $x_j \in \widehat{T} \backslash T$ that lies between $x_{r_i}$ and $x_{r_{i+1}}$. By definition of $\widehat{T}$, the neighbours of $x_j$ both lie on a common axis. In other words, we must have $r_{i+1} - r_i = 2$, and $x_{r_i}$ and $x_{r_{i+1}}$ belong to a common axis. Let now $\gamma \subseteq T$ be the path obtained by connecting the vertices of the form $x_{r_i}$ through the corresponding axes. Then $\gamma$ is a subtree of $T$ containing $x$ and $y$. It is convex, hence must contain the geodesic $\gamma_{x,y}^T$. This means we found a way to cover $\gamma$, and thus $\gamma_{x,y}^T$, with $D' \leq D$ axes. By definition of $M$ and $D$, we obtain
$$d_{\widehat{T}}(x,y) \geq M. \ \ (**)$$
By hypothesis, there is no axis of the form $h \cdot \gamma_s$ that covers a subgraph of $\gamma_{x,y}^T$ of diameter more than $K$. In particular then, one must use at least $d_T(x,y) / K$ such axes in order to cover $\gamma_{x,y}^T$ completely. This means $M \geq d_T(x,y) / K$. We conclude using $(**)$ that $d_{\widehat{T}}(x,y) \geq d_T(x,y) / K$, satisying $(*)$. This finishes the proof of the claim.
\bigskip

We now check that the hypothesis of the claim is satisfied. Suppose that no such constant $K$ exists. Then there is an axis $h \cdot \gamma_s$ such that the subtree $\gamma \coloneqq \gamma_g \cap h \cdot \gamma_s$ has diameter at least $2 \cdot ||g||+1$, where $||g||$ is the translation length of $g$ when acting on $T$. Since $\gamma$ has diameter at least $2 \cdot ||g||+1$, there is an edge $e \subseteq \gamma$ whose distance in $T$ to any of the two enpoints of $\gamma$ is at least $||g||$. Note that $e$ is a segment of $\gamma_g$, so $g \cdot e$ belongs to $\gamma_g$ as well. By definition, the distance in $T$ between $e$ and $g \cdot e$ is at most $||g||$, which means that $g \cdot e$ belongs to $\gamma$ as well. In particular, $g \cdot e$ belongs to $h \cdot \gamma_s$. Note that the action of $g$ respects the bipartite structure of $T$, and thus its translation length $||g||$ is an even number. Note on the other hand that the translation length of $h \cdot s \cdot h^{-1}$ when acting on $T$ is exactly $2$ (because its translation length when acting on $Y$ is $1$). Since $g \cdot e$ belongs to $h \cdot \gamma_s$, this means there is some constant $M$ such that $g \cdot e$ coincides with $h \cdot s^M \cdot h^{-1} \cdot e$ (actually, $M = \pm ||g||/2$). We get the equation
$$g \cdot e = h \cdot s^M \cdot h^{-1} \cdot e.$$
In particular, the element $g^{-1} \cdot h \cdot s^M \cdot h^{-1}$ stabilises $e$, hence must belong to $Z(A_{ab})$ by Lemma \ref{LemmaTTree}. We obtain $g = h \cdot s^M \cdot h^{-1} \cdot W$ for some $W \in Z(A_{ab})$. This is absurd by hypothesis.
\hfill\(\Box\)
\bigskip

\subsection{The syllabic length of powers of the Garside element}

We are now interested in the study of the elements of $A_{ab}$ that act elliptically on $\widehat{T}$, which have been described in Lemma \ref{LemmaEllipticOfHatT}. Our goal will be to give a linear lower bound on the syllabic length of powers of the Garside element (see Lemma \ref{LemmaLowerBoundForDelta}). The method is more algebraic, and we decide to briefly recall how one can obtain the Garside normal form of an element $g \in A_{ab}$ (see (\cite{mairesse2006growth}, Section 4) for a similar description).

\begin{algo} \label{Algo} Let $g \in A_{ab}$, and let $u \in F_{ab}$ be any word satisfying $u \equiv g$. Then one can obtain $\operatorname{Gars}(g)$ from $u$ in two steps:
\bigskip

\noindent \underline{Step 1:} If there is no occurence of a subword of the form $\Delta_x^{\pm1}$ in $u$, for some $x \in \{a, b\}$, or if all such occurences appear consecutively on the right-most part of $u$, go to Step 2. Otherwise, consider the left-most occurence of a $\Delta_x^{\pm1}$ subword in $u$, and write
$$u = v_1 \cdot \Delta_x^{\pm1} \cdot v_2,$$
for the appropriate subwords $v_1, v_2 \in F_{ab}$. Let
$$u' \coloneqq v_1 \cdot \widetilde{v_2} \cdot \Delta_x^{\pm1},$$
and note that $u' \equiv u \equiv g$. Then replace $u$ with $u'$, and proceed through Step 1 again.
\bigskip

\noindent \underline{Step 2:} At this point, we have a word $u$ that doesn't contain any subword of the form $\Delta_x^{\pm1}$, except potentially on its right-most part. This means $u$ takes the form
$$u = u_1 \cdots u_n \cdot W,$$
where each $u_i$ is an atom or the inverse of atom, and $W$ is a product of terms of the form $\Delta_a^{\pm1}$ and $\Delta_b^{\pm1}$. Moreover, for every $1 \leq i \leq n-1$, the last letter of $u_i$ and the first letter of $u_{i+1}$ either have opposite sign, or agree. If there is no negative letter (i.e. $a^{-1}$ or $b^{-1}$) in $u_1 \cdots u_n$, terminate the algorithm. Otherwise, the word $u_1 \cdots u_n$ contains at least one subword that is the inverse of an atom. Locate the left-most subword $u_i$ of this form. Without loss of generality, $u_i = (a^{-1}, b^{-1}; k)$ for some $1 \leq k < m$ (if $u_i$ starts with $b^{-1}$ instead, proceed symmetrically). Write
$$u = u_1 \cdots u_{i-1} \cdot (a^{-1}, b^{-1}; k) \cdot u_{i+1} \cdots u_n \cdot W,$$
and let
$$u' \coloneqq u_1 \cdots u_{i-1} \cdot (b,a;m-k) \cdot \widetilde{u_{i+1}} \cdots \widetilde{u_n} \cdot \widetilde{W} \cdot \Delta_x^{-1}$$
for some $x \in \{a, b\}$. One can check that $u' \equiv u \equiv g$. Replace $u$ by $u'$, and proceed through Step 2 again.
\end{algo}

\begin{ex} \label{Example} Let $m \coloneqq 3$, and let $u \coloneqq a b a^2 b^{-1} a^{-1} b a b a^2 b^4 a b$. We denote by $u_i$ the word obtained after the $i$-th Step of Algorithm \ref{Algo}. Then:
\begin{align*}
u_1 &= b a^{-1} b^{-1} a^2 b^3 \Delta_a \Delta_b \Delta_a \\
u_2 &= b^4 a^3 \Delta_b \Delta_a
\end{align*}
If we decompose the resulting word according to Definition \ref{DefiGarside}, we obtain
$$\operatorname{Gars}(u) = b \cdot b \cdot b \cdot ba \cdot a \cdot a \cdot \Delta_b \Delta_a = b^4 a^3 b a b a b a.$$
\end{ex}


\begin{lemma} \label{LemmaDeltanContainsADelta} Let $u \in F_{ab}$ and suppose that $u \equiv \Delta^n$ for some $n \neq 0$. Then $u$ contains a subword of the form $\Delta_x^{\pm 1}$ for some $x \in \{a, b\}$.
\end{lemma}

\noindent \textbf{Proof:} The proof uses the strategy of Algorithm \ref{Algo}. Suppose that $u$ does not contain any subword of the form $\Delta_x^{\pm 1}$ with $x \in \{a, b\}$. By definition, when giving $u$ as an imput, the first step of Algorithm \ref{Algo} is trivial. Starting with the second step of the algorithm, this means we can decompose $u$ in a product of atoms, inverses of atoms, and a power of the Garside element (see Algorithm \ref{Algo}):
$$u = u_1 \cdots u_k \cdot W$$
When applying the second step of the algorithm until the algorithm terminates, every atom $u_i$ yields an atom $u_i'$ that is either $u_i$ or $\widetilde{u_i}$, and every inverse of an atom $u_i$ yields an atom $u_i'$ that is either $u_i^*$ or $\widetilde{u_i}^*$, where $u_i^*$ is the unique atom such that $u_i = u_i^* \cdot \Delta_x^{-1}$ for some $x \in \{a, b\}$. Note that for every $1 \leq i \leq k$, $u_i$ is trivial if and only if $u_i'$ is trivial. We obtain the Garside normal form of $\Delta^n$:
$$\operatorname{Gars}(\Delta^n) = u_1' \cdots u_k' \cdot W',$$
for an appropriate $W'$. Recall that one can find trivial Garside normal forms for $\Delta^n$, such as $\operatorname{Gars}(\Delta^n) = \Delta_x^n$ for $x \in \{a, b\}$. By unicity of the atoms in the decomposition of $\operatorname{Gars}(\Delta^n)$, we obtain that all the $u_i'$ are trivial, and thus so are the $u_i$. In particular, $u = 1$, which is absurd.
\hfill\(\Box\)

\begin{lemma} \label{LemmaLowerBoundForDelta} For any $n \in \textbf{Z}$, $\ell_{\mathcal{S}}(\Delta^n) \geq (m-2) \cdot |n|$.
\end{lemma}

\noindent \textbf{Proof:} This is clear if $n = 0$. Since $\ell_{\mathcal{S}}(\Delta^n) = \ell_{\mathcal{S}}(\Delta^{-n})$, it is enough to prove that the result holds for $n > 0$. Let $u \in F_{ab}$ be any word representing $\Delta^n$. It is enough to show that
$$\ell_{\mathcal{S}}(u) \geq (m-2) \cdot n.$$
We now consider the string of words $u_0, u_1, u_2, \cdots u_{\lambda} \in F_{ab}$ defined by induction as follows. We first set $u_0 \coloneqq u$. By Lemma \ref{LemmaDeltanContainsADelta}, $u_0$ contains a subword of the form $\Delta_{x_1}^{\pm 1}$ for some $x_1 \in \{a, b\}$, so we can decompose $u_0$ as
$$u_0 = u_{0,1} \cdot \Delta_{x_1}^{\pm 1} \cdot u_{0,2},$$
for the appropriate words $u_{0,1}, u_{0,2} \in F_{ab}$. We then set
$$u_1 \coloneqq u_{0,1} \cdot \widetilde{u_{0,2}}.$$
Note that $u_1 \cdot \Delta_{x_1}^{\pm 1} \equiv u_0$. If $u_1$ is trivial, set $\lambda = 1$ and stop here. Otherwise, $u_1 \equiv \Delta^{n \pm 1}$ with $n \pm 1 \neq 0$, so we can apply Lemma \ref{LemmaDeltanContainsADelta} again and follow the same construction as above and obtain a word $u_2$ satisfying $u_2 \cdot \Delta_{x_2}^{\pm 1} = u_1$ for some $x_2 \in \{a, b\}$. As long as $u_i \not\equiv 1$, we continue to construct words $u_{i+1}$ in the fashion described above. The words obtained satisfy $u_{i+1} \cdot \Delta_{x_{i+1}}^{\pm 1} = u_i$ for some $x_{i+1} \in \{a, b\}$. Note that
$$u_i \ = \ \underbrace{ \overbrace{u_{i,1}}^{k_1 \text{ syl.}} \cdot \overbrace{\Delta_{x_i}^{\pm1}}^{m \text{ syl.}} \cdot \ \overbrace{u_{i,2}}^{k_2 \text{ syl.}} }_{\geq k_1 + k_2 + m - 2 \text{ syl.}}, \ \ \text{ and } \ \
u_{i+1} \ = \ \underbrace{u_{i,1} \cdot \widetilde{u_{i,2}}}_{\leq k_1 + k_2 \text{ syl.}},$$
so eventually
$$\ell_{\mathcal{S}}(u_i) \geq \ell_{\mathcal{S}}(u_{i+1}) + (m-2).$$
This means each word $u_{i+1}$ is syllabically shorter than $u_i$ by at least $(m-2)$ syllables. In particular, this process has to stop after a finite number $\lambda$ of steps. The final word, $u_{\lambda}$, satisfies
$$u_{\lambda} \cdot \prod\limits_{i=1}^{\lambda} \Delta_{x_i}^{\pm 1} \equiv u \equiv \Delta^n.$$
In particular then, $u_{\lambda}$ represents a power of $\Delta$, but does not contain any subword of the form $\Delta_x^{\pm 1}$ for some $x \in \{a, b\}$. By Lemma \ref{LemmaDeltanContainsADelta}, this means $u_{\lambda}$ is the trivial word. We obtain
$$\prod\limits_{i=1}^{\lambda} \Delta_{x_i}^{\pm 1} \equiv \Delta^n \Longrightarrow \lambda \geq n.$$
Trying to sum up the previous arguments, we have:
\\ (1) For $0 \leq i \leq \lambda -1$, each $u_{i+1}$ is syllabically shorter than $u_i$ by at least $(m-2)$ syllables. In particular, $u_{\lambda}$ is syllabically shorter than $u$ by at least $\lambda (m-2)$ syllables.
\\ (2) $u_{\lambda}$ is trivial.
\\ (3) $\lambda \geq n$.
\\Altogether, this gives a bound on the syllabic length of $u$:
$$\ell_{\mathcal{S}}(u) \overset{(1)}\geq \ell_{\mathcal{S}}(u_{\lambda}) + \lambda (m-2) \overset{(2)}= \lambda (m-2) \overset{(3)}\geq (m-2) \cdot n.$$
\hfill\(\Box\)

We are now able to prove the main Propositions:
\bigskip

\noindent \textbf{Proof of Proposition E.(2) and Proposition \ref{PropEquivToPropE}.(2):} We first recall that the two statements are equivalent, thanks to Lemma \ref{LemmaSyl}. Therefore we will only care on proving Proposition \ref{PropEquivToPropE}.(2). We divide the proof in four different cases. In all cases except the first one, we will give a linear lower bound of $\ell_{\mathcal{S}}(g^n)$ in terms of $n$. In all that follows, $h$ is an element of $A_{ab}$, $s \in \{a, b\}$ is a standard generator, and $W$ is an element of the centre $Z(A_{ab})$.
\bigskip

\noindent \underline{Case 1: $g = h \cdot s^k \cdot h^{-1}$.} Let $M \coloneqq \ell_{\mathcal{S}}(h) = \ell_{\mathcal{S}}(h^{-1})$. Then for any $n \in \textbf{Z}$, we have
$$\ell_{\mathcal{S}}(g^n) = \ell_{\mathcal{S}}(h \cdot s^{kn} \cdot h^{-1}) \leq \ell_{\mathcal{S}}(h) + \ell_{\mathcal{S}}(s^{kn}) + \ell_{\mathcal{S}}(h^{-1}) = M + 1 + M = 2M+1.$$

\noindent \underline{Case 2: $g = h \cdot s^k \cdot h^{-1} \cdot W$ with $W \neq 1$.} Then there is a $q \neq 0$ such that $g = h \cdot s^k \cdot h^{-1} \cdot \Delta^q$. Let $g_0 \coloneqq h \cdot s^k \cdot h^{-1}$, then $g^n = (g_0 \cdot \Delta^q)^n = g_0^n \cdot \Delta^{qn}$. On one hand we know by Case 1 that $\ell_{\mathcal{S}}(g_0^n)$ is uniformly bounded for all $n \geq 0$.  On the other hand, $\ell_{\mathcal{S}}(\Delta^{qn})$ grows linearly in $n$, by Lemma \ref{LemmaLowerBoundForDelta}. Putting these two facts together shows that $\ell_{\mathcal{S}}(g^n)$ grows linearly as well.
\bigskip

\noindent \underline{Case 3: $\exists N \neq 0: g^N \in Z(A_{ab})$.} By hypothesis, there is a $q \neq 0$ such that $g^N =\Delta^q$. By Lemma \ref{LemmaLowerBoundForDelta}, this means the quantity $\ell_{\mathcal{S}}(g^{Nn})$ grows linearly in $n$. In particular, the quantity $\ell_{\mathcal{S}}(g^{N \cdot \lfloor \frac{n}{N} \rfloor})$ grows linearly in $n$ as well (for a smaller constant). Note that the difference between $\ell_{\mathcal{S}}(g^n)$ and $\ell_{\mathcal{S}}(g^{N \cdot \lfloor \frac{n}{N} \rfloor})$ is uniformly bounded by the constant $L \coloneqq \max \{ \ell_{\mathcal{S}}(g^i) \ | \ i=0, \cdots, N-1 \}$. It follows that $\ell_{\mathcal{S}}(g^n)$ also grows linearly in $n$.
\bigskip

\noindent \underline{Case 4: We are in none of the previous cases.} Then by Lemma \ref{LemmaEllipticOfHatT}, $g$ acts loxodromically on $\widehat{T}$. In particular, the quantity $d_{\widehat{T}}(1_{\bullet}, g^n \cdot 1_{\bullet})$ grows linearly. We conclude with Lemma \ref{LemmaDistanceAndSyllabicLengthInHatT}.
\hfill\(\Box\)

\section{On the geometry of the action}

Let $A_{\Gamma}$ be a $2$-dimensional Artin group of rank at least $3$, and let $D_{\Gamma}$ be its modified Deligne complex. Our goal is to show that there exists a vertex $v \in D_{\Gamma}$, and an element $g \in A_{\Gamma}$ satisfying the two hypotheses of Theorem D. We have seen in Proposition E.(1) that a strong enough condition for $v$ to satisfy the first hypothesis of Theorem D is to require that its local group $G_v$ is a dihedral Artin group $A_{ab}$ with coefficient $3 \leq m_{ab} < \infty$. When such a vertex $v$ exists, it only remains to show that there exists an element $g \in A_{\Gamma}$ such that $A_{ab} \cap g A_{ab} g^{-1}$ is finite (which is equivalent to trivial because those groups are torsion-free). Our main geometric tool in order to find such an element is the following lemma:

\begin{lemma} \label{StabGeod} Let $G$ be a group acting by simplicial isomorphisms on a CAT(0) simplicial complex $X$ of dimension $2$. Let $v \in X$, $g \in G$ and denote by $G_p$ the stabiliser of a point $p \in X$. If the unique geodesic $\gamma$ between $v$ and $g v$ goes through a point with trivial stabiliser, then $G_v \cap G_{gv} = \{1\}$.
\end{lemma}

\noindent \textbf{Proof:} Any element of $G_v \cap G_{gv}$ fixes $v$ and $gv$, hence fixes (pointwise) the unique geodesic $\gamma$ between them. This means that $G_v \cap G_{gv} = G_{\gamma}$. Let $p \in \gamma$ be a point with trivial stabiliser. Then we have
$$G_v \cap G_{gv} = G_{\gamma} \subseteq G_p = \{1\}.$$
\hfill\(\Box\)

\noindent \textbf{Strategy:} The strategy of this section is led by the previous lemma. It is not hard to see that if $v \in D_{\Gamma}$ is a vertex with stabiliser $A_{ab}$, then the stabiliser of $gv$ for some $g \in A_{\Gamma}$ is exactly $g A_{ab} g^{-1}$. Suppose additionally that $v$ satisfies Theorem D.(1) (i.e. that $A_{ab}$ is large). Our goal will be to construct a geodesic between $v$ and some $gv$ that contains a point with trivial stabiliser. In the case of the modified Deligne complex, every point that lies in the interior of a triangle $T_{st}$ or an edge $e_s$ or $e_{st}$ has trivial stabiliser, hence it is enough to show that $\gamma$ goes through the interior of such a triangle or edge. In some cases, this will turn out to be quite difficult to prove. However, everything will be more manageable when working in some augmented version of the Deligne complex (see Definition \ref{DefAugDeligne}).
\bigskip

\noindent The next Proposition will give the structure of the different cases we will encounter:

\begin{prop} \label{PropCases} Let $A_{\Gamma}$ be a $2$-dimensional Artin group of rank at least $3$, and suppose that $\Gamma$ is connected and that $A_{\Gamma}$ is not a right angled Artin group. Then there exist three distinct generators $a,b,c \in S$ such that $m_{ab} \in \{3,4,\cdots \}$, $m_{ac} \in \{2,3,4,\cdots \}$, $m_{bc} \in \{2,3,4,\cdots,\infty\}$ and
$$\frac{1}{m_{ab}} + \frac{1}{m_{ac}} + \frac{1}{m_{bc}} \leq 1,$$
where $\frac{1}{\infty} \coloneqq 0$. Moreover, we are in exactly one of the following situation:
\\(1) There is a triplet $(a,b,c)$ as before that satisfies $m_{bc} < \infty$.
\\(2) There is no triplet $(a,b,c)$ as before with $m_{bc} < \infty$, but there is one that satisfies $m_{bc} = \infty$. Moreover, the graph $\Gamma^{bc}$ obtained from $\Gamma$ by adding an edge $e^{bc}$ with coefficient $6$ is such that $A_{\Gamma^{bc}}$ has dimension $2$.
\\(3) We are not in the first two situations, and $\Gamma$ contains a full cycle with coefficients $(2,2,2,n)$ for some $n \geq 3$.
\end{prop}

\noindent \textbf{Proof :} We begin by proving the first statement. Because $A_{\Gamma}$ is not right angled, there is an edge $e^{ab}$ in $\Gamma$ with coefficient $m_{ab} \in \{3,4,\cdots\}$. As $\Gamma$ is connected and has at least $3$ vertices, $e^{ab}$ has a neighboring edge in $\Gamma$, say $e^{ac}$, with coefficient $m_{ac} \in \{2,3,4,\cdots\}$. Since $A_{\Gamma}$ has dimension $2$, the last coefficient $m_{bc} \in \{2,3,4,\cdots,\infty\}$ satisfies:
$$\frac{1}{m_{ab}} + \frac{1}{m_{ac}} + \frac{1}{m_{bc}} \leq 1.$$
Let's now prove that we are in exactly one of the three cases. The three cases are exclusive by definition, so it is enough to show that if we are not in one of the first two situations, then we must be in the third. To prove this, pick a triplet of the form $m_{ab} \in \{3,4,\cdots \}$, $m_{ac} \in \{2,3,4,\cdots \}$, $m_{bc} = \infty$. By hypothesis, the graph $\Gamma^{bc}$ obtained from $\Gamma$ by adding an edge $e^{bc}$ with coefficient $6$ is such that $A_{\Gamma^{bc}}$ is not $2$-dimensional. This means that there is a generator $d \in S$ such that
$$\frac{1}{6} + \frac{1}{m_{bd}} + \frac{1}{m_{cd}} > 1.$$
This is only possible if $m_{bd} = m_{cd} = 2$. Notice that $m_{ad} = \infty$, otherwise the triplet $(a,b,d)$ would satisfy (1). This means that we have a full cycle $(e^{bd},e^{cd},e^{ac},e^{ab})$ with coefficients $(2,2,\geq 2,\geq 3)$ in $\Gamma$. If $m_{ac} = 2$, we are done. Suppose that $m_{ac} \geq 3$, and add an edge $e^{ad}$ of coefficient $6$. Since $A_{\Gamma^{ad}}$ is not $2$-dimensional by hypothesis and since
$$\frac{1}{6} + \frac{1}{m_{ab}} + \frac{1}{m_{bd}} \leq 1,$$
$$\frac{1}{6} + \frac{1}{m_{ac}} + \frac{1}{m_{cd}} \leq 1,$$
then there must be a fifth generator $e \in S$ such that
$$\frac{1}{6} + \frac{1}{m_{ae}} + \frac{1}{m_{de}} > 1.$$
For the same reasons as before, we have $m_{ae} = m_{de} = 2$ and $m_{ce} = \infty$. Hence there is a full cycle $(e^{ae},e^{de},e^{cd},e^{ac})$ with coefficients $(2,2,2,\geq 3)$ in $\Gamma$.
\hfill\(\Box\)
\bigskip

Recall that our goal in order to prove Theorem A is to apply Theorem D. For an irreducible $2$-dimensional Artin group $A_{\Gamma}$ of rank at least $3$, it turns out that the modified Deligne complex $D_{\Gamma}$ is exactly the space that we want to act on, at least in the first and third cases or Proposition \ref{PropCases}. Unfortunately, in the second case of Proposition \ref{PropCases}, the space $D_{\Gamma}$ is not fit to apply our main geometric tool that is Lemma \ref{StabGeod}. The reason, as will be seen later, is that we would like to have three generators $a,b,c \in S$ for which all the triangles $T_{ab}$, $T_{ba}$, $T_{bc}$, $T_{cb}$, $T_{ca}$,  and $T_{ac}$ belong to $D_{\Gamma}$. This is not the case when $m_{bc} = \infty$. However, notice that in the second case of Proposition \ref{PropCases}, the complex obtained from $D_{\Gamma}$ by adding the vertices of the form $g v_{bc}$ and their attached triangles $g T_{bc}$, $g T_{cb}$ is $2$-dimensional by hypothesis. This slightly bigger complex, as defined in the next definition, will be the one to look at when using Lemma \ref{StabGeod} and Theorem D in that case.

\begin{defi} \label{DefAugDeligne}  Let $A_{\Gamma}$ be a $2$-dimensional Artin group of rank at least $3$ with modified Deligne complex $D_{\Gamma}$ and fundamental domain $K_{\Gamma}$. Let $\Gamma^{st}$ be the same graph as $\Gamma$, except that we add an edge $e^{st}$ with coefficient $6$ between $s$ and $t$ if $m_{st} = \infty$. Consider now the $2$-dimensional complex $K_{\Gamma^{st}}$ obtained from Definition \ref{DefDeligne} for the group $A_{\Gamma^{st}}$. In other words,
$$K_{\Gamma^{st}} \coloneqq \begin{cases} K_{\Gamma} \text{ if } m_{st} < \infty \\
K_{\Gamma} \cup T_{st} \cup  T_{ts} \text{ if } m_{st} = \infty, \end{cases}$$
where the angle at $v_{st}$ in $T_{st}$ or $T_{ts}$ is set to be $\frac{\pi}{12}$ if $m_{st} = \infty$. We now want to realise $A_{\Gamma}$ as the fundamental group of a complex of groups over $K_{\Gamma^{st}}$. Doing so, we will give $K_{\Gamma^{st}}$ a structure of complex of groups, which may differ from the one coming from Definition \ref{DefDeligne}. When $m_{st} < \infty$, $K_{\Gamma^{st}} = K_{\Gamma}$, and we proceed as in Definition \ref{DefDeligne}: $K_{\Gamma^{st}}$ is simply the usual complex of groups associated with $A_{\Gamma}$. When $m_{st} = \infty$ however, $K_{\Gamma}$ is a strict subcomplex of $K_{\Gamma^{st}}$, and we carry the usual complex of groups associated with $A_{\Gamma}$ from $K_{\Gamma}$ onto the corresponding subcomplex of $K_{\Gamma^{st}}$. We still have to describe the local group at $v_{st}$ and the associated maps. We just set this local group to be the free group $F_{st}$ of rank $2$. The associated maps are the obvious morphisms that inject $\{1\}$, $\langle s \rangle$ and $\langle t \rangle$ into $F_{st}$. We call this structure of complex of group given to $K_{\Gamma^{st}}$ the \textbf{augmented} complex of groups associated with $A_{\Gamma}$ (relatively to $s$ and $t$).

If $m_{st} < \infty$, the augmented complex of groups associated with $A_{\Gamma}$ coincides with its usual complex of groups. If $m_{st} = \infty$, the augmented complex of groups associated with $A_{\Gamma}$ is the same as the one we would get if we took the usual complex of groups associated with $A_{\Gamma^{st}}$, but then replaced the local group $A_{st}$ at $v_{st}$, that is a dihedral Artin group with coefficient $6$, by the free group $F_{st}$. Note that in both cases, the $2$-dimensional complex under the augmented complex of groups associated with $A_{\Gamma}$ is $K_{\Gamma^{st}}$. The universal cover $D_{\Gamma}^{st}$ of this complex of groups is called the \textbf{augmented Deligne complex} of $A_{\Gamma}$ (relatively to $s$ and $t$). In particular, in the light of Definition \ref{DefDeligne}, we have
$$D_{\Gamma}^{st} \coloneqq \quotient{A_{\Gamma} \times K_{\Gamma^{st}}}{\sim},$$
where $(g,x) \sim (g',x') \Longleftrightarrow x = x'$ and $g^{-1} g'$ belongs to the local group of the smallest simplex of $K_{\Gamma^{st}}$ that contains $x$. The action of $A_{\Gamma}$ on itself induces an action of $A_{\Gamma}$ on $D_{\Gamma}^{st}$ by simplicial morphisms with strict fundamental domain $K_{\Gamma^{st}}$.
\end{defi}

\begin{figure}[H]
\centering
\includegraphics[scale=0.5]{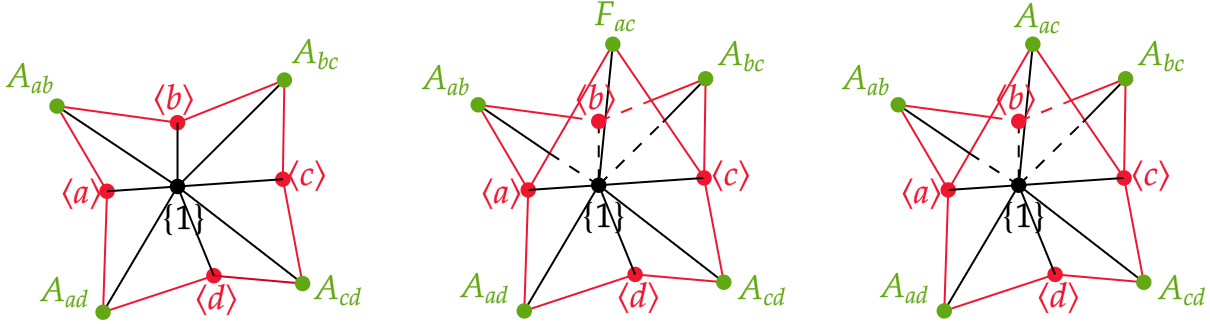}
\caption{Let $\Gamma$ be a square graph (as in Figure \ref{FigModDelCpx}) such that $m_{ac} = \infty$. Then are drawn:
\\ \underline{On the left:} The usual complex of groups associated with $A_{\Gamma}$.
\\ \underline{In the centre:} The augmented complex of groups associated with $A_{\Gamma}$ relatively to $a$ and $c$.
\\ \underline{On the right:} The usual complex of groups associated with $A_{\Gamma^{ac}}$.
\\Note that the first two complexes of groups share the same fundamental groups, and the last two complexes of groups share the same underlying $2$-dimensional complex.}
\label{FigureAugmentedComplex}
\end{figure}

\begin{rem} \label{RemDelAugmOk} (1) If $m_{st} <\infty$, the augmented Deligne complex $D_{\Gamma}^{st}$ and the modified Deligne complex $D_{\Gamma}$ agree.
\\(2) If $m_{st} = \infty$, then $D_{\Gamma}^{st}$ differs from $D_{\Gamma^{st}}$, as the fundamental groups of their associated complexes of groups are not the same: the former is $A_{\Gamma}$ while the latter is $A_{\Gamma^{st}}$. In particular, $D_{\Gamma}^{st}$ decomposes as a quotient of $A_{\Gamma} \times K_{\Gamma^{st}}$, while $D_{\Gamma^{st}}$ decomposes as a quotient of $A_{\Gamma^{st}} \times K_{\Gamma^{st}}$. Note however that the fundamental domains of these complexes are the same, as $2$-dimensional complexes.
\\(3) It is important to notice that if $a, b, s, t$ are four (non-necessarily all distinct) generators of $A_{\Gamma}$ satisfying $(a,b) \neq (s,t)$ and $m_{ab} < \infty$, then $Lk_{D_{\Gamma}}(v_{ab}) \cong Lk_{D_{\Gamma}^{st}}(v_{ab})$ (see Figure \ref{FigureAugmentedComplex} for instance). In particular, results such as Lemma \ref{LemmaSyl} or Proposition E.(1) also hold for $v_{ab}$ if we replace $D_{\Gamma}$ by $D_{\Gamma}^{st}$.
\end{rem}

The following lemma is a basic result but we decide to write it explicitly as it will be used several times in this section:

\begin{lemma} (\cite{bridson2013metric}, Chapter II.5) \label{CAT(0)sys} Let $X$ be a piecewise euclidean $2$-dimensional simply connected simplicial complex with finite number of shapes. Then $X$ is CAT(0) if and only if for every vertex $v \in X$ we have $\operatorname{sys}(Lk_X(v)) \geq 2 \pi$, where $\operatorname{sys}(Lk_X(v))$ is the length of the systole in $Lk_X(v)$.
\end{lemma}

\begin{lemma} \label{LemmaAugCAT(0)} Let $A_{\Gamma}$ be a $2$-dimensional Artin group of rank at least $3$, and suppose that we are in the second case of Proposition \ref{PropCases}. Then the augmented Deligne complex $D_{\Gamma}^{bc}$ of $A_{\Gamma}$ is CAT(0).
\end{lemma}

\noindent \textbf{Proof:} By hypothesis $A_{\Gamma^{bc}}$ has dimension $2$, hence its associated modified Deligne complex $D_{\Gamma^{bc}}$ is CAT(0) (Theorem \ref{XCat0}). We want to show that $D_{\Gamma}^{bc}$ is CAT(0). By Lemma \ref{CAT(0)sys}, and up to reducing to the fundamental domain, it is enough to show that every vertex $v \in K^{\Gamma^{bc}}$ satisfies
$$\operatorname{sys}(Lk_{D_{\Gamma}^{bc}}(v)) \geq 2 \pi. \ (*)$$
Notice that if $v \neq v_{bc}$ then
$$Lk_{D_{\Gamma}^{bc}}(v) \cong Lk_{D_{\Gamma^{bc}}}(v),$$
and thus $(*)$ follows from the fact that $D_{\Gamma^{bc}}$ is CAT(0), along with Lemma \ref{CAT(0)sys}.

If $v = v_{bc}$, then the local group at $v$ is the free group $F_{bc}$ by definition. We can do a similar analysis as the one done in Remark \ref{RemLinks}. This time, the maps of the development inject into the free group $F_{bc}$. Therefore, the link $Lk_{D_{\Gamma}^{bc}}(v_{bc})$ is isomorphic to the barycentric subdivision of the Bass-Serre tree above the segment of groups with local groups $\langle b \rangle$ and $\langle c \rangle$ on the vertices and $\{1\}$ on the edge. In particular, $Lk_{D_{\Gamma}^{bc}}(v_{bc})$ is simply-connected, i.e. has infinite systole.
\hfill\(\Box\)
\bigskip

We are now ready to prove the following lemma, that shows the existence of appropriate weakly malnormal subgroups of $A_{\Gamma}$, which is one of the requirements of Theorem D.

\begin{lemma} \label{LemmaGeod} Let $A_{\Gamma}$ be a $2$-dimensional Artin group of rank at least $3$, and suppose that $\Gamma$ is connected and that $A_{\Gamma}$ is not a right angled Artin group. Then there exists an Artin subgroup $A_{ab}$ with coefficient $3 \leq m_{ab} < \infty$ and an element $g \in A_{\Gamma}$ such that $A_{ab} \cap g A_{ab} g^{-1} = \{1\}$.
\end{lemma}

\noindent \textbf{Proof:} By Proposition \ref{PropCases}, we know that we either have three generators $a, b, c \in V(\Gamma)$ that satisfy exactly one of the following:
\\(1) $m_{ab}, m_{ac} \in \{3,4,\cdots\}$ and $m_{bc}  \in \{3,4,\cdots,\infty\}$;
\\(2) $m_{ac} = 2$, $m_{ab} \in \{3,4,\cdots\}$ and $m_{bc} \in \{5,6,\cdots,\infty\}$;
\\(3) $m_{ac} = 2$, $m_{ab} = m_{bc} = 4$.
\\Or we have four generators $a, b, c, d \in V(\Gamma)$ satisfying:
\\(4) The cycle $(e^{bc},e^{cd},e^{ad},e^{ab})$ is full in $\Gamma$ and has coefficients $(2,2,2,n)$ with $n \geq 3$.

\begin{figure}[H]
\centering
\includegraphics[scale=0.5]{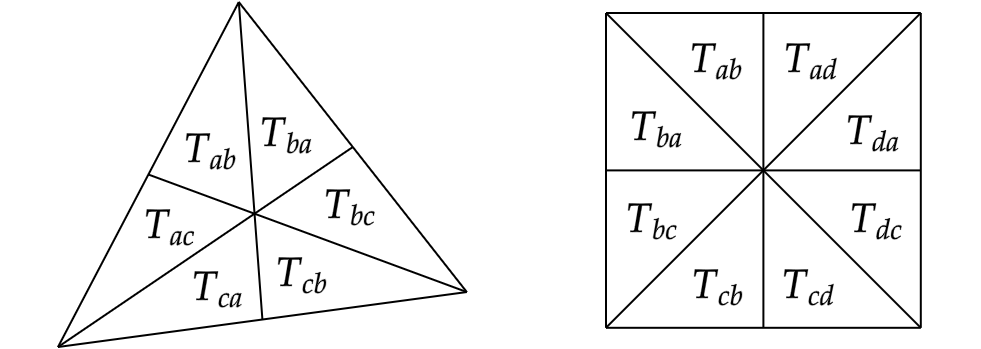}
\caption{$\Delta$ in the case (1), (2), (3) (on the left) and (4) (on the right).}
\label{Delta}
\end{figure}

\noindent Let $\Delta$ be the abstract complex (see Figure \ref{Delta}) defined by:
\\$\bullet$ In the situations (1), (2) and (3), $\Delta \coloneqq T_{ab} \cup T_{ba} \cup T_{bc} \cup T_{cb} \cup T_{ca} \cup T_{ac}$.
\\$\bullet$ In the situation (4), $\Delta \coloneqq T_{ab} \cup T_{ba} \cup T_{bc} \cup T_{cb} \cup T_{cd} \cup T_{dc} \cup T_{da} \cup T_{ad}$.
\\Note that in either case, the points in the interior of $\Delta$ have trivial stabilisers. Also note that we don't have to look at the augmented Deligne complex in the situations (1) and (2) if $m_{bc} < \infty$, and neither do we in the situations (3) and (4). However in those cases $m_{bc} < \infty$, and hence $D_{\Gamma}^{bc} = D_{\Gamma}$, so it will just be convenient to write $D_{\Gamma}^{bc}$ to cover all cases.

Let now $P$ be an abstract complex defined by $P \coloneqq \quotient{P_0}{\sim}$, where $P_0$ is defined depending on the situations given in the beginning of the proof by:
\\(1) $P_0 \coloneqq \Delta \sqcup (c \Delta)$.
\\(2) $P_0 \coloneqq \Delta \sqcup (c \Delta) \sqcup (cb \Delta) \sqcup (cbc \Delta)$.
\\(3) $P_0 \coloneqq \Delta \sqcup (c \Delta) \sqcup (cb \Delta) \sqcup (cbc \Delta) \sqcup (cba \Delta) \sqcup (cbca \Delta) \sqcup (cbcab \Delta) \sqcup (cbcabc \Delta)$.
\\(4) $P_0 \coloneqq \Delta \sqcup (c \Delta) \sqcup (d \Delta) \sqcup (cd \Delta)$.

\begin{figure}[H]
\centering
\includegraphics[scale=0.5]{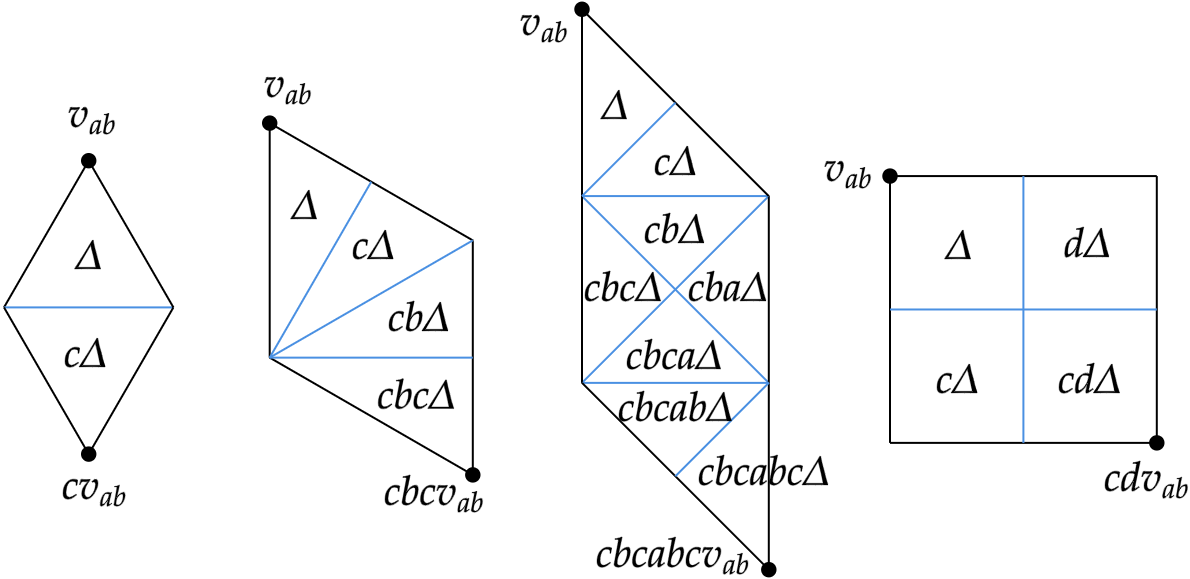}
\caption{Polygon $P$ in the four different cases, from left to right.}
\label{P}
\end{figure}

\noindent And where $\sim$ corresponds to the gluing shown in Figure \ref{P}, i.e. $P$ is obtained from $P_0$ by gluing the different copies of $\Delta$ along some of their edges (drawn in blue in Figure \ref{P}).

\noindent A priori, we can't be sure that there are no additional gluings happening in $D_{\Gamma}^{bc}$, so we don't want to look at $P$ as a subcomplex of $D_{\Gamma}^{bc}$, but we want instead look at $P$ through the natural map $f: P \rightarrow D_{\Gamma}^{bc}$ that maps $P$ to $D_{\Gamma}^{bc}$.  
\bigskip

\noindent \underline{\textbf{Claim 1:}} $P$ is isometrically embedded in $D_{\Gamma}^{bc}$.

\noindent \textbf{Proof:} In the light of (\cite{charney2000tits}, Lemma 1.4), it is enough to show that for every $p \in P$, the induced map $f_p: Lk_P(p) \rightarrow Lk_{D_{\Gamma}^{bc}}(p)$ is $\pi$-distance preserving, i.e. that
$$\forall x, y \in Lk_P(p), \ d_{Lk_P(p)}(x,y) \geq \pi \Rightarrow d_{Lk_{D_{\Gamma}^{bc}}(p)}(f_p(x),f_p(y)) \geq \pi.$$
There are two different situations:
\\$\bullet$ If $p \in P$ is in the orbit of $v_{\emptyset}$, then $Lk_{D_{\Gamma}^{bc}}(p)$ is just the augmented defining graph $\Gamma^{bc}$ with the appropriate metric (see Definition \ref{DefAugDeligne} and Remark \ref{RemLinks}). Notice that, any edge $e^{st} = e_{s,st} \star e_{t,st}$ from $s$ to $t$ in $\Gamma^{bc}$ has length
$$\ell(e^{st}) = 2 \cdot \angle_{\emptyset}(v_s,v_{st}) = \pi - \frac{\pi}{m_{st}} \geq \frac{\pi}{2},$$
according to the metric on $Lk_{D_{\Gamma}^{bc}}(p)$. Since $Lk_P(p)$ is simply the full cycle in $\Gamma^{bc}$ corresponding to the triangle $(e^{ab},e^{ac},e^{bc})$ (in the situations (1), (2) and (3)) or to the square $(e^{bc}, e^{cd}, e^{ad}, e^{ab})$ (in the situation (4)), we can apply (\cite{charney2000tits}, Lemma 1.6) and conclude that the map $f_p: Lk_P(p) \hookrightarrow Lk_{D_{\Gamma}^{bc}}(p)$ is $\pi$-distance preserving.
\\$\bullet$ If $p \in P$ is not in the orbit of $v_{\emptyset}$, then it is not hard to see from Remark \ref{RemLinks} that every full cycle in $Lk_P(p)$ has length exactly $2 \pi$. In particular, the map $f_p$ must be $\pi$-preserving, otherwise we would be able to build an isometrically embedded cycle in $Lk_{D_{\Gamma}^{bc}}(p)$ of length strictly less than $2 \pi$, contradicting the CAT(0)-ness of $D_{\Gamma}^{bc}$ (Theorem \ref{XCat0}, Lemma \ref{LemmaAugCAT(0)} and Lemma \ref{CAT(0)sys}).

We can now use (\cite{charney2000tits}, Lemma 1.4) and conclude that $P$ is isometrically embedded in $X$. In particular, geodesics in $P$ project to geodesics in $X$ through $f$.
\bigskip

\noindent \underline{\textbf{Claim 2:}} $A_{ab} \cap g A_{ab} g^{-1} = \{1\}$ for some $g \in A_{\Gamma}$.

\noindent \textbf{Proof:} Notice that $P$ is CAT(0) by Lemma \ref{CAT(0)sys}. In particular, it is uniquely geodesic. Let $\gamma$ be the geodesic in $P$ defined depending on the situations given in the beginning of the proof by:
\\(1) $\gamma$ is the geodesic going from $v_{ab}$ to $c v_{ab}$.
\\(2) $\gamma$ is the geodesic going from $v_{ab}$ to $cbc v_{ab}$.
\\(3) $\gamma$ is the geodesic going from $v_{ab}$ to $cbcabc v_{ab}$.
\\(4) $\gamma$ is the geodesic going from $v_{ab}$ to $cd v_{ab}$.
\\Note that $\gamma$ is also geodesic in $D_{\Gamma}^{bc}$, by the previous claim. Thanks to Lemma \ref{StabGeod}, it is enough to show that $\gamma$ goes through the interior of some $g_0 \Delta$ contained in $P$. Consider either of the four situations and suppose that it is not the case. Then in particular $\gamma$ would be contained in the $1$-skeleton of $P$. It is not hard to check, since we know every angle in $P$ by construction, that there must be a vertex $v$ in $\gamma$ that satisfies $\angle_v^P(\gamma) < \pi$. This is not possible, as $\gamma$ is a geodesic and $P$ is CAT(0).
\hfill\(\Box\)
\bigskip

We have worked through everything that was required in order to use our main criterion, that is Theorem D. We can now prove our main Theorem:
\bigskip

\begin{theorem*} Every irreducible $2$-dimensional Artin group of rank at least $3$ is acylindrically hyperbolic.
\end{theorem*}

\noindent \textbf{Proof of Theorem A:} Let $A_{\Gamma}$ be an irreducible $2$-dimensional Artin group of rank at least $3$. We can assume that $\Gamma$ is connected, as otherwise $A_{\Gamma}$ splits as a free product $A_{\Gamma_1} * A_{\Gamma_2}$ of infinite groups hence is acylindrically hyperbolic. We can also assume that $A_{\Gamma}$ is not a right angled Artin group, as every irreducible right angled Artin group that is not cyclic is acylindrically hyperbolic (\cite{osin2016acylindrically}, 8.(d)).

Let $a, b, c \in V(\Gamma)$ be the three generators obtained in the proof of Lemma \ref{LemmaGeod}, and consider the action of $A_{\Gamma}$ on its augmented Deligne complex $D_{\Gamma}^{bc}$. Note that the latter is CAT(0) by Lemma \ref{LemmaAugCAT(0)} and Lemma \ref{XCat0}. Since $m_{ab} \geq 3$, we know from Proposition E.(1) and Remark \ref{RemDelAugmOk} that the orbits of $A_{ab}$ on $Lk_{D_{\Gamma}^{bc}}(v_{ab})$ are unbounded. Moreover, we know from Lemma \ref{LemmaGeod} that there exists an element $g \in A_{\Gamma}$ such that $A_{ab} \cap g A_{ab} g^{-1} = \{1\}$. Therefore, we can apply Theorem D and conclude that $A_{\Gamma}$ is either virtually cyclic or acylindrically hyperbolic. That $A_{\Gamma}$ is not virtually cyclic is clear because it contains dihedral Artin subgroups that contain $\mathbf{Z}^2$.
\hfill\(\Box\)

\nocite{*}

\newcommand{\etalchar}[1]{$^{#1}$}

\bigskip
  \footnotesize

 \textsc{Department of mathematics, Heriot-Watt University, Riccarton, Edinburgh EH14 4AS} \par\nopagebreak
 \text{E-mail address:}
 \texttt{\href{mailto: ncv1@hw.ac.uk}{ncv1@hw.ac.uk}}

\end{document}